\documentclass[12pt]{article}
\usepackage{epsfig}
\usepackage{amsmath,amsfonts}
\title{JSJ-Decompositions of Coxeter Groups over Virtually Abelian Splittings}
\author{M. Mihalik}
\newtheorem{theorem}{Theorem}
\newtheorem{proposition}[theorem]{Proposition}
\newtheorem{lemma}[theorem]{Lemma}
\newtheorem{corollary}[theorem]{Corollary}
\newcounter{remarknum}
\newenvironment{remark}{\addvspace{12pt}\refstepcounter{remarknum}
\noindent{\bf Remark \arabic{remarknum}.}}{\par\addvspace{12pt}}
\newenvironment{proof}{\addvspace{10pt}\noindent{\bf Proof:}}{
$\Box$\par\addvspace{10pt}}
\newcounter{examplenum}
\newenvironment{example}{\addvspace{12pt}\refstepcounter{examplenum}
\noindent{\bf Example \arabic{examplenum}.}}{\par\addvspace{12pt}}
\date{December 3, 2009 
}
\begin{document}
\maketitle

\begin{abstract}
The idea of  ``JSJ-decompositions" for 3-manifolds began with work of Waldhausen and was developed later through work of Jaco, Shalen and Johansen. It was shown that there is a finite collection of 2-sided, incompressible tori that separate a closed irreducible 3-manifold into pieces with strong topological structure. 

Sela introduced JSJ-decompositions for groups, an idea that has flourished in a variety of directions. The general idea is to consider a certain class $\mathcal G$ of groups and splittings of groups in $\mathcal G$ by groups in another class $\mathcal C$. E.g. Rips and Sela considered splittings of finitely presented groups by infinite cyclic groups. For an arbitrary group $G$ in $\mathcal G$, the goal is to produce a ``unique"  graph of groups decomposition $\Psi$ of $G$ with edge groups in $\mathcal C$ so that $\Psi$ reveals all reduced graph of groups decompositions of $G$ with edge groups in $\mathcal C$. More specifically, if $H$ is a vertex group of $\Psi$ then either there is no $\mathcal C$-group that splits both $G$ and $H$, or $H$ has a special ``surface group-like" structure. Vertex groups of the second type are standardly called orbifold groups.

For a finitely generated Coxeter system $(W,S)$, we produce a reduced JSJ-decomposition $\Psi$ for splittings of $W$ over virtually abelian subgroups. We show $\Psi$ is unique (up to conjugate vertex groups) and each vertex and edge group is generated by a subset of $S$ (and so $\Psi$ is ``visual"). The construction of $\Psi$ is algorithmic. If $V\subset S$ generates an orbifold (vertex) group then $V=T\cup M$, where $\langle M\rangle$ is virtually abelian,  $\langle T\rangle$ is virtually a closed surface group or virtually free and $\langle V\rangle =\langle M\rangle \times \langle T\rangle$ . 
\end{abstract}


\section{Introduction}\label{Intro}

The theory of JSJ-decompositions has its origins in the work of Waldhausen \cite{Wald} on characteristic submanifolds of a 3-manifold and later work of Jaco-Shalen \cite{JS} and Johansen \cite{J}. For a closed, irreducible, oriented 3-manifold there is a finite collection of embedded 2-sided incompressible tori that separate the manifold into pieces, each of which is a Seifert fibered space or an atoroidal and acylindrical space. This gives a graph of groups decomposition of the fundamental group with edge groups free abelian of rank 2. 

Sela \cite{S} introduced the notion of JSJ-decomposition for a general class of groups and showed that word hyperbolic groups have JSJ-decompositions over infinite cyclic splittings. Rips-Sela \cite{RS} generalize this to finitely presented groups. Scott-Swarup \cite{SS} consider splittings corresponding to virtually polycyclic groups with restrictions on Hirsch length and extend these results to virtually abelian groups of bounded rank. Dunwoody-Sageev \cite{DS} and then Fujiwara-Papasoglu \cite {FP} gave JSJ-decompositions for finitely presented groups over slender splittings.

A group is {\it slender} if all of its subgroups are finitely
generated. The class of slender groups is contained in the class
of {\it small} groups which are defined in terms of actions on
trees. If a group contains a non-abelian free group it is not
small. Coxeter groups containing no non-abelian free group are in
fact virtually abelian and decompose in a special way amenable to
our results (see theorem \ref{T1}).

In analogy with the 1-ended assumptions of
Rips-Sela \cite{RS}, and following Dunwoody-Sageev \cite{DS}
directly, we define the class of minimal virtually abelian
splitting subgroups of $W$. If a Coxeter group splits over a minimal subgroup that contains no non-abelian free group, then the splitting subgroup is in fact virtually abelian. Hence for our purposes, there is no difference between (minimal) splittings over small, slender or virtually abelian groups and we only consider splittings of Coxeter groups over virtually
abelian subgroups.

If $(W,S)$ is a finitely generated
Coxeter system, a graph of groups decomposition $\Psi$ of $W$ is
{\it visual} if each edge and vertex group is generated by some
subset of $S$ and the bonding maps are inclusions. The main
theorem of \cite{MT} states that for any graph of groups
decomposition $\Lambda$ of $W$ there is a visual decomposition
$\Psi$ such that each vertex (respectively edge) group of $\Psi$
is conjugate to a subgroup of a vertex (respectively edge) group
of $\Lambda$. This result is used extensively in this paper and
the basics of visual decompositions are reviewed in section
\ref{Basic}.

Our construction of a JSJ-decomposition of a finitely generated
Coxeter group $W$, with virtually abelian edge groups
begins with a graph of groups decomposition $\Psi_1$ of $W$ with
edge groups that are minimal virtually abelian splitting
subgroups of $W$ and such that $\Psi_1$ is a maximal such
decomposition without edge groups that are ``crossing splitters".
We also show that $\Psi_1$ is unique (up to conjugate vertex
groups) and visual. We call $\Psi_1$ a {\it level 1
JSJ-decomposition of $W$ with virtually abelian edge groups}. If
$V\subset S$ and $\langle V\rangle$ is a vertex group of
$\Psi_1$, then $\langle V\rangle$ may contain virtually
abelian subgroups that split $\langle V\rangle$ and $W$ non-trivially. Minimal splitting subgroups of this type are not necessarily minimal
virtually abelian splitting subgroups of $W$. It is shown that
any graph of groups decomposition of $\langle V\rangle$ with such edge groups is compatible with
$\Psi_1$ and a maximal such decomposition that avoids crossing splitters, is
unique (up to conjugate vertex groups) and visual. Replacing all
such vertex groups of $\Psi_1$ by such graph of groups
decompositions gives $\Psi_2$, a {\it level 2 JSJ-decomposition of
$W$ with virtually abelian edge groups}. Continuing, we eventually
have a visual graph of groups decomposition $\Psi$ such that if
$\langle V\rangle$ is a vertex group of $\Psi$, then the only
virtually abelian splitting subgroups of $\langle V\rangle$ (that also split $W$) are crossing.
We call $\Psi$ a {\it JSJ-decomposition of $W$ with respect to
virtually abelian splittings} and show that $\Psi$ is unique (up
to conjugate vertex groups). If $\Lambda$ and $\Phi$ are graph of
groups decompositions of a group $W$, then say {\it the decomposition of $\Lambda$ induced by $\Phi$ is compatible with $\Lambda$} if
for each vertex group $V$ of $\Lambda$, the decomposition of $V$
induced by the action of $V$ on the Bass-Serre tree for $\Phi$ is
compatible with $\Lambda$. In sections \ref{JSJ} and \ref{Orbi},
we prove results that imply our main theorem:

\begin{theorem} \label{main} Suppose $(W,S)$ is a finitely
generated Coxeter system and $\Psi$ the reduced JSJ-decomposition
of $W$ with virtually abelian edge groups. Then:
\begin{enumerate}

\item $\Psi$ is visual, unique up to conjugate vertex groups, and algorithmically defined.

\item If $\Phi$ is a graph of groups decomposition of $W$ with
virtually abelian edge groups then the decomposition of $\Psi$
induced by $\Phi$ is compatible with $\Psi$.

\item If both $W$ and a vertex group $\langle V\rangle$ of $\Psi$ ($V\subset S$)
split nontrivially over a virtually abelian subgroup of $\langle V\rangle$, then
$\langle V\rangle$ decomposes as $\langle T\rangle\times \langle
M\rangle$ where $T\cup M=V$, $M$ generates a virtually abelian
group and the presentation diagram of $T$ is either a loop of
length $\geq 4$ (in which case $T$ generates a group that is
virtually a closed surface group) or the presentation diagram of
$T$ is a disjoint union of vertices and simple paths (in which
case $T$ generates a virtually free group with graph of groups
decomposition such that each vertex group is either $\mathbb Z_2$
or finite dihedral and each edge group is either trivial or
$\mathbb Z_2$).
\end{enumerate}
\end{theorem}

Vertex groups of the type described by part 3 of theorem \ref{main} are called {\it orbifold} vertex groups. If $H$ is a non-orbifold vertex group of our JSJ-decomposition $\Psi$, then $H$ does not split non-trivially
over a virtually abelian subgroup that also splits $W$ non-trivially. In particular, if $\Phi$ is another graph of
groups decomposition of $W$ with virtually abelian edge groups,
then $H$ is a subgroup of a conjugate of a vertex group of $\Phi$.

The construction of our JSJ-decompositions with virtually abelian
edge groups is algorithmic. Given the presentation diagram
$\Gamma$ of a Coxeter system $(W,S)$, theorem \ref{T1} allows us
to determine the subsets of $S$ that generate virtually abelian
subgroups. Those that separate $\Gamma$, algebraically split $W$.
A result in section \ref{Def} allows us to easily decide which of
these are the visual minimal virtually abelian splitting
subgroups of $(W,S)$ at all stages of the construction of the
$i^{th}$-level JSJ-decompositions. It is equally easy to
``visually" determine which of these splitting subgroups are
crossing and hence build JSJ-decompositions. Our main result
distinguishes the two types of vertex groups of our JSJ-decompositions. Those with no crossing subgroups are
indecomposable with respect to virtually abelian splittings of $W$ and
those with crossing subgroups which are traditionally called 
orbifold vertex groups.

If $\Phi$ is a graph of groups decomposition for a group $G$, $V$ is a vertex of $\Phi$ with incident edge $E$, and the groups of $V$ and $E$ agree, then one can collapse the decomposition $\Phi$ across the edge $E$ to obtain a smaller (more reduced) decomposition of $G$. As a simple example consider the splitting $A\ast _C E\ast _ED$ that collapses to $A\ast_CD$. While edges such as $E$ seem to contribute somewhat artificial splittings to $\Phi$, this type of edge is important for the JSJ decompositions produced by Fujiwara and Papasoglu. The advantage of the unreduced splitting $A\ast _C E\ast _ED$ is that it exhibits the decomposition $\langle A\cup E\rangle \ast_E D$ whereas $A\ast _C D$ does not. Our decompositions are reduced. The connection between our decompositions and the Fujiwara/Papasoglu decompositions is that their decompositions collapse to ours. 

In \cite{GL}, Guirardel-Levitt, develop the idea that a JSJ decomposition should be a deformation space satisfying a universal property. A deformation space (introduced by Forester in \cite{Fo}) is a collection of $G$-trees, which in fact is a contractible complex. In the correct setting, the deformation spaces of Guirardel-Levitt contain the trees constructed in \cite{DS}, \cite{FP} and \cite{RS}, as well as the ones constructed in this paper. 
For Coxeter groups the trees for the Fujiwara/Papasoglu decomposition and our decomposition are in the same deformation space. 

\section{Basic Facts and Background Results}\label{Basic}

A thorough discussion of graphs of groups decompositions of
Coxeter groups is given in \cite{MT}. We briefly discuss the
aspects of this theory necessary to this paper. Every Coxeter
group has a set of order 2 generators and so there is no
non-trivial map of a Coxeter group to $\mathbb Z$. In particular,
no Coxeter group is an HNN extension of any sort and any graph of
groups decomposition of a Coxeter group has graph a tree. Hence
the decompositions of Coxeter groups are a straightforward
generalization of amalgamated product decompositions. For a graph
of groups decomposition $\Lambda$ of a group $G$, the Bass-Serre
tree $T$ for $\Lambda$ has vertices (respectively edges) the
cosets $wV$ where $w\in G$ and $V$ is a vertex (respectively
edge) group of $\Lambda$.  There is a left action of $G$ on $T$ and an element $g$ of $G$ stabilizes the
coset $wV$ iff $g\in wVw^{-1}$. If $V$ is a vertex of $\Lambda$ with
vertex group $\Lambda (V)$, and $\Phi$ a graph of groups
decomposition of $\Lambda (V)$, then $\Phi$ is {\it compatible}
with $\Lambda$ if for each edge $E$ of $\Lambda$ incident to $V$,
$\Lambda (E)$ is contained in a $\Lambda (V)$-conjugate of a vertex group of
$\Phi$. In this case $V$ can be replaced by $\Phi$ to produce a
finer graph of groups decomposition of $G$. A graph of groups
decomposition $\Lambda$ is {\it reduced} if no edge between
distinct vertices has edge group the same as an end point vertex
group. If a graph of groups is not reduced, we may collapse a
vertex group across an edge, where the edge group is the same as
the endpoint vertex group, giving a smaller graph of groups
decomposition of the original group.

\begin{lemma} \label{Conj}
(See \cite{MT}) Suppose $\Lambda$ is a reduced graph of groups decomposition of a
group $G$, the underlying graph for $\Lambda$ is a tree, $U$ and $V$ are vertices of $\Lambda$, and
$g\Lambda(U)g^{-1}\subset \Lambda (V)$ for some $g\in G$ then $U=V$ and $g\in \Lambda (U)$. $\square$
\end{lemma}

The following result is straightforward.
\begin{lemma} \label{Compat}
Suppose $\Psi$ is a graph of groups decomposition of a group $G$,
$V$ is a finitely generated vertex group of $\Psi$, and $\cal E$
is a collection of subgroups of $V$ such that for any $K\in \cal E
$, $V$ splits non-trivially and $\Psi$-compatibly over $K$. If
$\Lambda$ is a graph of groups decomposition of $V$ with edge
groups in $\cal E$, then $\Lambda$ is compatible with $\Psi$.
$\square$
\end{lemma}

The next result easily follows from the combinatorics of
group actions on trees or more practically from the exactness of
the Mayer-Vietoris sequence for a pair of groups.

\begin{lemma} \label{NoZ}
Suppose a group $G$ splits as $A\ast _CB$. If there
is no non-trivial homomorphism from $G$ or $C$ to $\mathbb Z$,
then there is no non-trivial homomorphism from $A$ or $B$ to
$\mathbb Z$. In particular, if $\Lambda$ is a graph of groups
decomposition of a Coxeter group and no edge group of $\Lambda$
maps non-trivially to $\mathbb Z$, then no vertex group of
$\Lambda$ maps non-trivially to $\mathbb Z$. $\square$
\end{lemma}

We take a {\it Coxeter presentation} to be given as
$$P=\langle S:(st)^{m(s,t)}\ (s,t\in S,\,m(s,t)< \infty)\rangle$$
where $m:S^2\to\lbrace 1,2,\ldots,\infty\rbrace$ is such that
$m(s,t)=1$ iff $s=t$, and $m(s,t)=m(t,s)$.  In the group with this
presentation, the elements of $S$ represent distinct elements of
order 2 and a product $st$ of generators has order $m(s,t)$.  A
{\it Coxeter group} $W$ is a group having a Coxeter presentation
and a {\it Coxeter system} $(W,S)$ is a Coxeter group $W$ with
generating subset $S$ corresponding to the generators in a Coxeter
presentation of $W$. When the order of the product of a pair of
generators is infinite there will be no defining relator for that
pair of generators and we will say that the generators are {\it
unrelated}. Our basic reference for Coxeter groups is Bourbaki
\cite{Bourbaki}. A {\it special} or {\it visual} subgroup for a
Coxeter system $(W,S)$, is a subgroup of $W$ generated by a
subset of $S$. If $W'$ is the visual subgroup generated by
$S'\subseteq S$ in a Coxeter system $(W,S)$, then $(W',S')$ is
also a Coxeter system. More specifically the following result
(see \cite{Bourbaki}) is fundamental to the study of Coxeter
groups.

\begin{proposition}\label{flat}
Suppose $(W,S)$ is a Coxeter system and $P=\langle S:
(st)^{m(s,t)}$ for $m(s,t)<\infty\rangle$ (where $m:S^2\to
\{1,2,\ldots ,\infty\}$) is a Coxeter presentation for $W$. If
$A\subset S$, then $(\langle A\rangle, A)$ is a Coxeter system
with Coxeter presentation $\langle A:(st)^{m'(s,t)}$ for
$m'(s,t)<\infty\rangle$ (where $m'=m\vert _{A^2}$). $\square$
\end{proposition}

Given a group $G$ and a generating set $S$, an {\it $S$-geodesic}
for $g\in G$ is a shortest word in $S\cup S^{-1}$ such that the
product of the letters of this word is $g$. The number of letters
in an $S$-geodesic for $g$ is the {\it $S$-length of $g$}.  An
important combinatorial fact about geodesics for a Coxeter system
is called the ``deletion condition".

\begin{proposition}\label{dc} {\bf (The Deletion Condition)} Suppose
$(W,S)$ is a Coxeter system and $w=a_1\cdots a_n$ for $a_i\in S$.
If $a_1\cdots a_n$ is not geodesic then there are indices $i<j$
in $\{1,2,\ldots ,n\}$ such that $w=a_1\cdots a_{i-1}a_{i+1}\cdots
a_{j-1}a_{j+1}\cdots a_n$. I.e. $a_i$ and $a_j$ can be deleted.
$\square$
\end{proposition}

The information given by a Coxeter presentation may be
conveniently expressed in the form of a labeled graph. We define
the {\it presentation diagram} of the system $(W,S)$ to be the
labeled graph $\Gamma(W,S)$ with vertex set $S$, and an
(undirected) edge labeled $m(s,t)$ between distinct vertices $s$
and $t$ when $m(s,t)<\infty$. The connected components of the
presentation diagram $\Gamma(W,S)$ correspond to visual subgroups
which are the factors in a free product decomposition of $W$. In
contrast, a {\it Coxeter graph} has vertex set $S$ and labeled
edges when $m(s,t)\neq 2$. The components of a Coxeter graph
corresponding to direct product factors of $W$. By proposition
\ref{flat}, a presentation diagram of a visual subgroup of $W$
generated by a subset $S'\subseteq S$ is the induced subgraph of
$\Gamma(W,S)$ with vertex set $S'$.

Suppose $\Gamma(W,S)=\Gamma_1\cup\Gamma_2$ is a union of induced
subgraphs and let $\Gamma_0=\Gamma_1\cap\Gamma_2$ (so vertices and
edges of $\Gamma(W,S)$ are in $\Gamma_1$ or $\Gamma_2$ or both,
and $\Gamma_0$ is the induced subgraph consisting of the vertices
and edges in both).  Equivalently, suppose $\Gamma_0$ is an
induced subgraph with $\Gamma(W,S)-\Gamma_0$ having at least two
components, $\Gamma_1$ is $\Gamma_0$ together with some of these
components and $\Gamma_2$ is $\Gamma_0$ together with the other
components. We say in this case that $\Gamma_0$ {\it separates}
$\Gamma(W,S)$ (separates it into at least two components).  Then
it is evident from the Coxeter presentation that $W$ is an
amalgamated product of visual subgroups corresponding to
$\Gamma_1$ and $\Gamma_2$ over the visual subgroup corresponding
to $\Gamma_0$. Amalgamated product decompositions with visual
factors and visual amalgamated subgroup are easily seen in the
presentation diagram and we call such an amalgamated product a
{\it visual} splitting of $W$.

  We say that $\Psi$ is a {\it visual graph of
groups decomposition} of $W$ (for a given Coxeter system $(W,S)$),
if each vertex and edge group of $\Psi$ is visual for $(W,S)$, the
injections of each edge group into its endpoint vertex groups are
given simply by inclusion, and the fundamental group of $\Psi$ is
isomorphic to $W$ by the homomorphism induced by the inclusion
map of vertex groups into $W$.  A sequence of compatible visual
splittings of $W$ will result in such a decomposition. In
\cite{MT}, we study general graph of groups decompositions of
Coxeter groups and how these are related to visual graph of groups
decompositions. The main result of \cite{MT} shows that an
arbitrary graph of groups decomposition of a Coxeter group can be
refined (in a certain sense) to a visual graph of groups
decomposition.

\begin{theorem}\label{MT1}
Suppose $(W,S)$ is a Coxeter system and $W$ is a subgroup of the
fundamental group of a graph of groups $\Lambda$. Then $W$ has a
visual graph of groups decomposition $\Psi$ where each vertex
group of $\Psi$ is a subgroup of a conjugate of a vertex group of
$\Lambda$, and each edge group of $\Psi$ is a subgroup of a
conjugate of an edge group of $\Lambda$.  Moreover, $\Psi$ can be
taken so that each visual subgroup of $W$ that is a subgroup of a
conjugate of a vertex group of $\Lambda$ is a subgroup of a
vertex group of $\Psi$. $\square$
\end{theorem}

The following three results are important technical facts proved
in \cite{MT}.

\begin{lemma}\label{visgog}
Suppose $(W,S)$ is a Coxeter system.  A graph of groups $\Psi$
with graph a tree, where each vertex group and edge group is a
visual subgroup for $(W,S)$ and each edge map is given by
inclusion, is a visual graph of groups decomposition of $W$ iff
each edge in the presentation diagram of $W$ is an edge in the
presentation diagram of a vertex group and, for each generator
$s\in S$, the set of vertices and edges with groups containing
$s$ is a nonempty subtree in $\Psi$. $\square$
\end{lemma}

If $\Psi$ is a visual graph of groups decomposition for the
Coxeter system $(W,S)$, it is convenient to label the vertices of $\Psi$ by the subsets of $S$ that generate the
corresponding vertex groups. So if $Q\subset S$ is a
vertex of $\Psi$, then $\Psi (Q)=\langle Q\rangle$. The same convention is not used for edge groups as distinct edges may have the same edge group.

Separation properties of various subsets of a presentation diagram
$\Gamma (W,S)$ for a Coxeter system $(W,S)$ are frequently
analyzed in this paper. If $A$ and $B$ are subsets of $\Gamma$ we
say $A$ {\it separates} $B$ {\it in} $\Gamma$ if there are points
$b_1$ and $b_2$ of $B-A$ such that any path in $\Gamma$ from
$b_1$ to $b_2$ intersects $A$.

\begin{lemma} \label{LC1}
Suppose $(W,S)$ is a Coxeter system, $S$ is finite, $\Psi$ is a
visual graph of groups decomposition of $W$, $E'$ is an edge of
$\Psi$ and $\Psi(E')=\langle E\rangle$ for $E\subset S$. If $\{x,y\}\subset S-E$, and $x\in X$ and $y\in Y$ for
$X$ and $Y$ vertices of $\Psi$ on opposite sides of $E'$, then $E$ separates $x$ and $y$ in $\Gamma (W,S)$.
$\square$
\end{lemma}

\begin{lemma}\label{LC2}
Suppose $(W,S)$ is a finitely generated Coxeter system, and
$\Psi$ is a visual graph of groups decomposition of $W$. If $C$ is
a complete subset of the presentation diagram $\Gamma (W,S)$, then
there is a vertex $V$ of $\Psi$ such that $C\subset V$. $\square$
\end{lemma}

In this paper we require more than just the statement of theorem
\ref{MT1}. The technique to produce a visual decomposition is easy
to describe and useful to the constructions in this paper. Under
the hypothesis of theorem \ref{MT1} let $T$ be the Bass-Serre tree
for $\Lambda$.  In the proof of theorem \ref{MT1}, it is shown
that $W$ has a visual graph of groups decomposition with graph
$T$ and vertex (respectively edge) group at $gV$ generated by the
subset of $S$ that stabilizes $gV$. Since $S$ is finite, this
visual graph of groups reduces to a finite reduced graph of groups
decomposition of $W$ satisfying the conclusion of theorem \ref{MT1}. In this paper we make repeated use of this
construction and refer to it as {\it the visual graph of groups
given by the construction for theorem} \ref{MT1}.

The next lemma follows from a result of Kilmoyer (see section 4
of \cite{MT}).

\begin{lemma} \label{L3}
Suppose $(W,S)$ is a Coxeter system, $I$, $J\subset S$, and $d$
is a minimal length double coset representative in $\langle
I\rangle w\langle J\rangle$. Then $\langle I\rangle\cap d\langle
J\rangle d^{-1}=\langle K\rangle$ for $K=I\cap (dJd^{-1})$ and,
$d^{-1}\langle K\rangle d=\langle J\rangle \cap (d^{-1}\langle
I\rangle d)=\langle K'\rangle$ for $K'=J\cap d^{-1}Id=d^{-1}Kd$.
In particular, if $w=idj$ for $i\in \langle I\rangle$ and $j\in
\langle J\rangle$ then $\langle I\rangle \cap w\langle J \rangle
w^{-1}=i\langle K\rangle i^{-1}$ and $\langle J\rangle \cap
w^{-1}\langle I\rangle w=j^{-1}\langle K'\rangle j$. $\square$
\end{lemma}

\section{Decomposing Coxeter Groups With No Non-Abelian Free Subgroup}\label{PT1}

The Euclidean simplex reflection groups and irreducible
finite Coxeter groups are catalogued in Coxeter's book \cite{Cox}.
Euclidean simplex groups are infinite and virtually abelian.
In his thesis \cite{Krammer}, D. Krammer classifies free abelian
subgroups of Coxeter groups. If a Coxeter groups contains no
non-abelian free group, the Coxeter group is virtually free abelian and
decomposes as a direct product of finite Coxeter groups and
Euclidean simplex groups in a special way - a result well-known
to experts. We include an elementary proof of this result, using
the following lemma (which is a direct consequence of theorems
7.2.2, 7.3.1, 12.1.19 and exercise 12.1.14 of J. Ratcliffe's book
\cite{Rat}).

\begin{lemma}\label{L2}
Consider the collection of Coxeter systems $(W,S)$ such that
\begin{enumerate}
\item $\Gamma (W,S)$ is complete,

\item $W$ does not decompose as $\langle A\rangle \times \langle
B\rangle$ for $A$ and $B$ non-trivial subsets of $S$, and

\item for each $s\in S$, $\langle S-\{s\}\rangle$ is either finite or
Euclidean.
\end{enumerate}
Then $W$ is finite, Euclidean or contains a free subgroup of rank
2. $\square$
\end{lemma}

\begin{theorem}\label{T1}
Suppose $(W,S)$ is a finitely generated Coxeter system and $W$
does not contain a non-abelian free group. Then $S$ is the
disjoint union of sets which commute with one another and each
generates a finite group or an Euclidean simplex group.
\end{theorem}

\begin{proof} 
If $a,b\in S$ such that $m(a,b)=\infty$, let
$C=S-\{a,b\}$. Then $W=\langle \{a\}\cup C\rangle \ast _{\langle
C\rangle }\langle \{b\}\cup C\rangle$. The index of $\langle
C\rangle$ in both $\langle \{a\}\cup C\rangle$ and $\langle
\{b\}\cup C\rangle$ is 2, since, $W$ contains no free group of
rank 2. Hence $\langle C\rangle$ is normal in $\langle \{a\}\cup
C\rangle$. For each $c\in C$, $aca\in \langle C\rangle$. Any
geodesic in $\langle C\rangle$ uses only letters in $C$, and by
the deletion condition, $aca=c$. Thus $a$ commutes with $C$ as
does $b$. The group $\langle a,b\rangle$ splits off as a direct
factor of $W$. Hence we may assume $\Gamma (W,S)$ is complete
(and infinite). Note that condition 1) of the lemma is satisfied.

Suppose $(W,S)$ is a counterexample to the theorem, with $\vert
S\vert$ small as possible. Note that $\vert S\vert\geq 3$, and
$W$ does not visually decompose as a non-trivial direct product,
so that condition 2) of the lemma is satisfied. If $s\in S$, we
have $\langle S-\{s\}\rangle$ is a visual product of a finite
group and Euclidean simplex groups. Choose $a$ and $b$ distinct
elements of $S$. Write $\langle S-\{a\}\rangle =\langle
F_a\rangle \times \langle E_1\rangle\times \cdots \times \langle
E_p\rangle$ (so $F_a\cup (\cup_{i=1} ^p E_i)=S-\{a\}$) where
$\langle F_a\rangle$ is finite and each $\langle E_i\rangle$ is
Euclidean. Similarly write $\langle S-\{b\}\rangle =\langle
F_b\rangle \times \langle K_1\rangle\times \cdots \times \langle
K_q\rangle$. If $b\in E_i$ for some $i$, then assume $i=1$. If
$p>1$, then $E_2\subset K_j$ for some $j$. For any proper subset
$K$ of $K_j$, $\langle K\rangle$ is finite, and so $E_2=K_j$. But
then $ E_2$ commutes with $S-E_2$ which is impossible. Hence
$\{p,q\}\subset \{0,1\}$.

If $F_a\ne \emptyset$ and $E_1\ne\emptyset$, we may choose $b\in
F_a$, so that $\{a,b\}\cap E_1=\emptyset$. This implies $E_1=K_1$
and $E_1$ commutes with $S-E_1$ $(=F_a\cup \{a\})$, which is
impossible. Hence either $F_a$ or $E_1$ is empty. Since $a$ is
arbitrary in $S$, condition 3) of the previous lemma is satisfied.
\end{proof}


\section{Minimal Virtually Abelian Splitting Subgroups}\label{Def}

If $\Lambda$ is a graph of groups decomposition of a group $W$
and $G$ is a vertex group of $\Lambda$, then a virtually abelian
subgroup $A$ of $G$ is a {\it minimal virtually abelian splitting
subgroup} for $(\Lambda ,G)$ if $G$ splits non-trivially and
compatibly with $\Lambda$ over $A$, and there is no virtually
abelian subgroup $B$ of $W$ such that $G$ splits non-trivially
and compatibly with $\Lambda$ over $B$, and $B\cap A$ has
infinite index in $A$ and finite index in $B$.

For a Coxeter system $(W,S)$, $\Psi$ a visual graph of groups
decomposition for $(W,S)$ and $G$ a vertex group of $\Psi$, let
${\cal C}(\Psi, G)$ be the set of virtually abelian 
subgroups of $G$ that split $G$ non-trivially and $\Psi$-compatibly. 
Let $M(\Psi,G)$ be the set of minimal virtually abelian
splitting subgroups for $(\Psi ,G)$.

If $\Psi$ is the trivial graph of groups decomposition for
$(W,S)$ (with one vertex), then define ${\cal C}(W,S)\equiv {\cal
C} (\Psi ,W)$ and $M(W,S)\equiv M(\Psi, W)$.
Observe:
\begin{enumerate}

\item If a vertex group $G$ of $\Psi$ has more than 1-end, then each
member of $M(\Psi ,G)$ is a finite group.

\item For a given finitely generated Coxeter group $W$, the ranks of
the virtually abelian subgroups of $W$ are bounded.

\item If $A\subset S$ and $\langle A\rangle \in M(\Psi ,G)$ then $A$
satisfies the conclusion of theorem \ref{T1}. Hence $A$ is the
disjoint union of sets that commute with one another such that (at most) one
generates a finite group and each other generates an Euclidean
simplex group.
\end{enumerate}

If $(W,S)$ is a Coxeter system and $A\subset S$ is such that
$\langle A\rangle$ is virtually abelian, define $E(A)$ to be the
set of generators of the Euclidean factors of the visual direct
product decomposition of $\langle A\rangle$ given by theorem
\ref{T1}.

If $A$ is a minimal virtually abelian splitting subgroup of a Coxeter
group, then by theorem \ref{MT1}, $A$ contains a subgroup of
finite index which is isomorphic to a Coxeter group. Hence there
is no non-trivial homomorphism of $A$ to $\mathbb Z$. Lemma
\ref{NoZ} implies the next result.

\begin{lemma}
Suppose $W$ is a finitely generated Coxeter group and $\Lambda$
is a graph of groups decomposition of $W$ with minimal virtually
abelian splitting subgroups as edge groups. Then $\Lambda$ is a tree and no vertex group
of $\Lambda$ maps non-trivially to $\mathbb Z$. $\square$
\end{lemma}

If $\Psi$ is a visual graph of groups decomposition for the
Coxeter system $(W,S)$, $V(\subset S)$ is a vertex of $\Psi$ and
$\Psi _1$ is a visual decomposition for $(\langle V\rangle, V)$,
then $\Psi _1$ is {\it visually compatible} with $\Psi$ if for
each edge group $\langle E\rangle$ ($E\subset S$) of an edge of $\Psi$ incident to $V$, $E$
is a subset of $K$ for some vertex $K(\subset V)$ of $\Psi _1$. In particular, if $A\subset V$ is such that $\langle A\rangle$ is an edge group of $\Psi_1$, then we say {\it $\langle A\rangle$ splits $V$ visually and compatibly with $\Psi$}. 
\newpage

\begin{lemma} \label{LM11}
Suppose $(W,S)$ is a finitely generated Coxeter system, $\Psi$ is
a reduced visual graph of groups decomposition for $(W,S)$,
$V\subset S$ is a vertex of $\Psi$, and $\Lambda $ is a reduced
graph of groups decomposition of $\langle V\rangle$ such that
each edge group of $\Lambda$ is in $M(\Psi ,\langle V\rangle)$.
Let $\Psi '$ be the reduced visual decomposition for $\Lambda$
given by the construction for theorem \ref{MT1}, then $\Psi '$ is
visually compatible with $\Psi$, each edge group of $\Psi '$ is
in $M(\Psi ,\langle V\rangle)$ and if $U\subset V$ such that
$\langle U\rangle$ is an edge group of $\Psi'$, then $U$
separates $V$ in $\Gamma (W,S)$.
\end{lemma}

\begin{proof} 
By lemma \ref{Compat}, $\Lambda$ is
compatible with $\Psi$. So if $E'$ is an edge of $\Psi
$ incident to $V$ and $E\subset S$ is such that $\langle E\rangle = \Psi(E')$, then $\langle E\rangle$ is a
subgroup of a $\langle V\rangle $-conjugate of a vertex group of
$\Lambda$. Hence $E$ stabilizes a vertex of the Bass-Serre tree
for $\Lambda$. But then by the construction for theorem \ref{MT1},
$E$ is a subset of a vertex group of $\Psi'$ and $\Psi'$ is
visually compatible with $\Psi$. Each edge group of $\Psi'$ is
conjugate to a subgroup of an edge group of $\Lambda$ and so must
be in $M(\Psi,\langle V\rangle)$. Suppose $A, B\subset V$ are
distinct vertices of $\Psi'$ incident to the edge $ U$ of
$\Psi'$. There exists $a\in A-U$ and $b\in B-U$. If $\Psi''$ is
the visual graph of groups decomposition obtained from $\Psi$ by
replacing the vertex $ V$ by $\Psi'$, then applying lemma
\ref{LC1} to $\Psi ''$ shows $ U$ separates $a$ and $b$ in
$\Gamma$. 
\end{proof}

The next lemma is a direct consequence of V.
Deodhar's results in \cite{Deo}.

\begin{lemma}\label{L4}
Suppose $(W,S)$ is a Coxeter system, $A\subset S$, $\langle
A\rangle$ is infinite and there is no non-trivial $F\subset A$
such that $\langle F\rangle$ is finite and $\langle A\rangle
=\langle A-F\rangle \times \langle F\rangle$. If $w\in W$ such
that $w\langle A\rangle w^{-1}\subset \langle B\rangle$ for
$B\subset S$, then $uAu^{-1}=A\subset B$ for $u$ the minimal
length double coset representative of $\langle B\rangle w\langle
A\rangle$. In particular, if $w\langle A\rangle w^{-1}=\langle
B\rangle$, then $A=B$. $\square$
\end{lemma}

Observe that if $(W,S)$ is a Coxeter system, and $W=\langle
F\rangle \times \langle G\rangle=\langle H\rangle \times \langle
I\rangle$ for $F\cup G=S=H\cup I$. Then $W=\langle F\rangle
\times \langle H-F\rangle\times \langle G\cap I\rangle$.

\begin{proposition} \label{P5}
Suppose $(W,S)$ is a Coxeter system, $A\subset S$ and $\langle
A\rangle=\langle B\rangle \times \langle C\rangle$ where $B\cup
C=A$, and $C$ is the (unique) largest such subset of $A$ such that
$\langle C\rangle $ is finite. If $w\langle A\rangle
w^{-1}\subset \langle H\rangle$ for $w\in W$ and $H\subset S$,
then $B\subset H$. In particular, if $\langle A\rangle$ and
$\langle H\rangle$ are virtually abelian then $E(A)\subset E(H)$.
\end{proposition}

\begin{proof} 
Let $d$ be a minimal length double coset
representative of $\langle H\rangle w\langle A\rangle$. So
$d\langle A\rangle d^{-1}\subset \langle H\rangle$ and $d\langle
A\rangle d^{-1}=d\langle A\rangle d^{-1}\cap \langle H\rangle
=\langle K\rangle$ for $K=dAd^{-1}\cap H$. Hence $K=dAd^{-1}$
($d^{-1}Kd\subset A$ and generates $\langle A\rangle$). Apply
lemma \ref{L4}. 
\end{proof}

\begin{figure}[ht]
\begin{center}
\epsfig{file=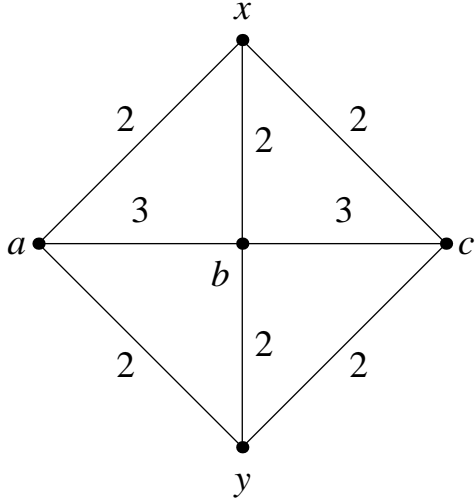}
\caption{$\Gamma(W,S)$ \label{jsjfig1}}
\end{center}
\end{figure}

\begin{example}\label{E1}
Consider the Coxeter system $(W,S)$ with presentation diagram shown in figure \ref{jsjfig1}. Clearly, $\langle
x,y,b\rangle\in M(W,S)$ and the element $bc$ conjugates
$\{x,y,b\}$ to $\{x,y,c\}$. So, $\langle x,y,c\rangle \in M(W,S)$.
Hence a visual subgroup in $M(W,S)$ need not separate the
presentation diagram $\Gamma (W,S)$.
\end{example}

\begin{figure}[ht]
\begin{center}
\epsfig{file=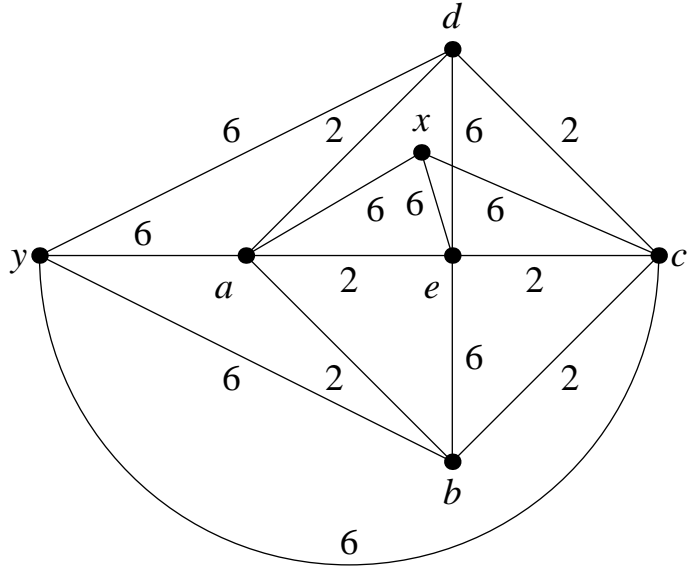}
\caption{$\Gamma(W,S)$ \label{jsjfig2}}
\end{center}
\end{figure}

\begin{example}\label{E2}
The group with presentation diagram shown in figure \ref{jsjfig2},
splits non-trivially and visually over both $\langle
a,b,c,d\rangle\in {\cal C}(W,S)$ and $\langle a,c,e\rangle$. As
$W$ is 1-ended and $\langle a,c,e\rangle$ is 2-ended, $\langle
a,c,e\rangle \in M(W,S)$. Hence $\langle a,b,c,d\rangle\not \in
M(W,S)$. Furthermore, no $A\subset \{a,b,c,d\}$ is such that
$\langle A\rangle \in M(W,S)$.
\end{example}

\begin{proposition}\label{P6}
Suppose $(W,S)$ is a Coxeter system, $\Psi$ is a visual graph of
groups decomposition for $(W,S)$ and $V\subset S$ a vertex of
$\Psi$. If $A\subset V$ and $\langle A\rangle\in {\cal C}(\Psi,
\langle V\rangle)$, then there exists $B\subset V$ such that $B$
separates $V$ in $\Gamma (W,S)$, $\langle B\rangle \in M(\Psi,
\langle V\rangle)$ and $E(B)\subset E(A)$.
\end{proposition}

\begin{proof} 
If $\langle A\rangle \in M(\Psi,\langle V\rangle)$,
let $G=\langle A\rangle$. Otherwise, let $G$ be a minimum rank
element of ${\cal C}(\Psi, \langle V\rangle )$ such that $G\cap
\langle A\rangle$ has infinite index in $\langle A\rangle$ and
finite index in $G$. In any case, $ G\in M(\Psi, \langle
V\rangle)$. By Lemma \ref{LM11}, there exists $B\subset V$ and
$v\in \langle V\rangle$ such that $B$ separates $V$ in $\Gamma$,
$v\langle B\rangle v^{-1}\subset G$ and $\langle B\rangle \in
M(\Psi,\langle V\rangle$). Hence $Rank(\langle B\rangle)=Rank
(G)\leq Rank(\langle A\rangle)$ and $\langle A\rangle \cap
v\langle B\rangle v^{-1}$ has finite index in $v\langle B\rangle
v^{-1}$.

By lemma \ref{L3}, $\langle A\rangle \cap v\langle B\rangle
v^{-1}= a\langle K\rangle a^{-1}$, for some $a\in \langle
A\rangle$ and $K\subset A$. Hence $E(K)\subset E(A)$. By lemma
\ref{L3}, $\langle B\rangle \cap v^{-1}\langle A\rangle v=b\langle
K'\rangle b^{-1}$ for some $b\in \langle B\rangle$ and $K'\subset
B$ where $K$ and $K'$ are conjugate. Hence by lemma \ref{L4}
$E(K)=E(K')$. As $b\langle K'\rangle b^{-1}$ has finite index in
$\langle B\rangle$, $E(K')=E(B)$. Hence $E(B)\subset E(A)$.
\end{proof}

We can now easily recognize $V$-separating visual subgroups in either
${\cal C}(\Psi, \langle V\rangle)$ or $M(\Psi, \langle V\rangle)$.

\begin{corollary}\label{C7}
Suppose $(W,S)$ is a Coxeter system, $\Psi$ is a reduced visual
graph of groups decomposition for $(W,S)$, $V\subset S$ is a
vertex group of $\Psi$, $A\subset V$ such that $\langle A\rangle$ splits $\langle V\rangle$ non-trivially and visually compatible with $\Psi$.
Then $\langle A\rangle \in M(\Psi, \langle V\rangle)$ iff
$\langle A\rangle $ is virtually abelian, and there is no
$B\subset V$ such that $\langle B\rangle$ is virtually abelian,
$\langle B\rangle$ splits $\langle V\rangle$ visually compatible with $\Psi$, and $E(B)$ is a proper subset of
$E(A)$. $\square$
\end{corollary}

For $(W,S)$ a Coxeter system, $\Psi$ a reduced visual
graph of groups decomposition for $(W,S)$, $V\subset S$ such that
$V$ is a vertex of $\Psi$, and $A\subset V$ such that $A$ separates $V$ in $\Gamma (W,S)$ and $ \langle A\rangle
\in M(\Psi, \langle V\rangle)$, we have not ruled out the possibility that there is $x\in E(A)$ such that $A-\{x\}$ separates $V$ in $\Gamma$. It may be that $x\in E$  for $\langle E\rangle$ the group of an edge 
incident to $V$ in $\Psi$ and that in $\Gamma$, $A-\{x\}$ separates $x$ from some other point of $E$, so that the visual splitting of $W$ over $\langle A-\{x\}\rangle$ is not compatible with $\Psi$ (and so the minimality of $\langle A\rangle$ is not violated). In our main applications, $\Psi$ will be an $n^{th}$-stage JSJ-decomposition and we will show there is no $B\subset V$ such that $\langle B\rangle$ is virtually abelian and $B$ separates $E$ in $\Gamma$ for $\langle E\rangle$ the group of an edge of $\Psi$ incident to $V$. But, until that point of the paper is reached, we add a restriction to the statements of some of our results in order to deal with this contingency.  

\begin{lemma} \label{L8}
Suppose $(W,S)$ is a Coxeter system, $\Psi$ is a reduced visual
graph of groups decomposition for $(W,S)$, $V\subset S$ is a vertex of $\Psi$, $A\subset V$ separates $V$ in $\Gamma (W,S)$, $\langle A\rangle$ is virtually abelian, and there is no $x\in E(A)$ such that $A-\{x\}$ separates $V$ in $\Gamma (W,S)$. If $K$ is a component of $\Gamma
-A$ which intersects $V$ non-trivially, then for each $a\in E(A)$
there is an edge from $K$ to $a$.
\end{lemma}

\begin{proof} 
Otherwise, $A-\{a\}$ separates $V$ in
$\Gamma$. 
\end{proof}

\begin{lemma} \label{L9}
Suppose $\Gamma$ is a graph with vertex set $S$, and $A\subset S$ separates $\{b_1,b_2\}$ in $\Gamma$. If $b\in S$ is adjacent to both $b_1$ and
$b_2$ in $\Gamma$, then $b\in A$. $\square$
\end{lemma}

\begin{lemma} \label{L10}
Suppose $(W,S)$ is a Coxeter system, $\Psi$ is a reduced visual
graph of groups decomposition for $(W,S)$, $V\subset S$ is a
vertex of $\Psi$, $A\subset S$ such that $\langle A\rangle$ is virtually abelian, and $B\subset
V$ such that $ \langle B\rangle$ is virtually abelian, $B$
separates $V$ in $\Gamma (W,S)$ and there is no $x\in E(B)$ such that $B -\{x\}$ separates $V$ in $\Gamma(W,S)$. If $A$ separates $B$ in $\Gamma$
then $B$ separates $A$ in $\Gamma$. In particular, (by lemma
\ref{L9}) $A-\{a_1,a_2\}=B-\{b_1,b_2\}$ for $a_1,a_2$ unrelated
elements of $A$ and $b_1,b_2$ unrelated elements of $B$. (So
$Rank(A)=Rank(B)$.) Furthermore, if $A\subset V$, $\langle A\rangle \in {\cal C}(\Psi , \langle V\rangle)$, $\langle B\rangle \in M(\Psi,\langle V\rangle)$, and for each $D\subset V$ such that $\langle D\rangle \in M(\Psi , \langle V\rangle)$ there is no $x\in E(D)$ such that $D-\{x\}$ separates $V$ in $\Gamma(W,S)$  then
$\langle A\rangle\in M(\Psi,\langle V\rangle)$.
\end{lemma}

\begin{proof} 
Write $B=\{b_1,b_2\}\cup M$, where $A$
separates $b_1$ and $b_2$ in $\Gamma$ and $\langle M\rangle$
commutes with $\{b_1,b_2\}$. By lemma \ref{L9}, $M\subset A$ and so
$A\cap B=M$. If $B$ does not separate $A$, then $A-M\subset K$
for $K$ a component of $\Gamma -B$.
For $i\in \{1,2\}$, let $K_i$ be distinct components of $\Gamma
-B$ containing $t_i\in V$. Assume $K_1\ne K$. Then $t_1\not \in
A$. We may assume $t_1$ is in a component of $\Gamma -A$ not
containing $b_2$. By hypothesis,  $\{b_1\}\cup M$ does not separate $V$ in $\Gamma$.
Choose a shortest path from $t_1$ to $t_2$ avoiding $\{b_1\}\cup
M$. Then this path passes through $b_2$ (exactly once). As $A$
separates $b_2$ and $t_1$, we let $s$ be the first vertex of $A$
in our path. Then $s\in A-M\subset K$ and we have connected
$t_1\in K_1$ to $A-M\subset K$ in $\Gamma -B$, which is
nonsense. By lemma \ref{L9}, $A-\{a_1,a_2\}=B-\{b_1,b_2\}\equiv M$.

It remains to show that if $A\subset V$, $\langle V\rangle \in {\cal C}(\Psi,\langle V\rangle)$ and $\langle B\rangle \in M(\Psi,\langle V\rangle)$ then $\langle A\rangle
\in M(\Psi, \langle V\rangle)$. Otherwise, proposition \ref{P6}
implies there is $D\subset V$ such that $D$ separates $V$ in
$\Gamma$, $\langle D\rangle \in M(\Psi, \langle V\rangle)$,
$E(D)\subset E(A)(= \{a_1,a_2\}\cup E(M))$ and $Rank (\langle D\rangle )<Rank(\langle
A\rangle)$. If $\{a_1,a_2\}\subset D$, then as $B$ separates
$a_1$ and $a_2$, the argument above shows
$D-\{a_1,a_2\}=B-\{b_1',b_2'\}$. But then $Rank (\langle
D\rangle)=Rank (\langle B\rangle)=Rank (\langle A\rangle)$ which
is nonsense. If $a_1\not\in D$, then $E(D)\subset M$. This is
impossible as $\langle B\rangle \in M(\Psi, \langle V\rangle)$
and $\langle B\rangle=\langle b_1,b_2\rangle \times \langle
M\rangle$. Similarly for $a_2$.
\end{proof}

If $\langle A\rangle $ and $\langle B\rangle$ are elements of
$M(\Psi, \langle V\rangle)$ and $A$ and $B$ separate one another
in $\Gamma$, we say $\langle A\rangle $ and $\langle B\rangle$
{\it cross} or are {\it crossing} in $M(\Psi, \langle V\rangle)$.

\begin{proposition} \label{P11}
Suppose $(W,S)$ is a Coxeter system, $\Psi$ is a reduced visual
graph of groups decomposition for $(W,S)$, $V\subset S$ is a
vertex of $\Psi$, $A\subset S$ such that $\langle A\rangle$ is virtually abelian, and $B\subset
V$ such that $ \langle B\rangle$ is virtually abelian, $B$
separates $V$ in $\Gamma (W,S)$ and there is no $x\in E(B)$ such that $B -\{x\}$ separates $V$ in $\Gamma(W,S)$. If $A$ separates $B$ in $\Gamma (W,S)$ then
$\Gamma (W,S)-B$ has exactly 2-components which intersect $V$
non-trivially.
\end{proposition}

\begin{proof} 
By lemma \ref{L10}, we may assume
$A=\{a_1,a_2\}\cup M$ and $B=\{b_1,b_2\}\cup M$, where $a_1$ and
$a_2$ are unrelated and separated by $B$, $b_1$ and $b_2$ are
unrelated and separated by $A$, and $M$ commutes with $\{a_1,a_2,
b_1,b_2\}$. For $i\in \{1,2\}$, let $K_i$ be the component of
$\Gamma -B$ containing $a_i$.

If $K$ is a component of $\Gamma -B$ which intersects $V$
non-trivially, then by lemma \ref{L8}, there is an edge connecting
$b_i$ to $K$ for $i\in \{1,2\}$. If additionally, $K\ne K_i$ for
$i\in \{1,2\}$, then $K\cap A=\emptyset$, and so $K\cup
\{b_1,b_2\}$ is a connected subset of $\Gamma -A$. This is
impossible since $A$ separates $b_1$ and $b_2$ in $\Gamma$.
\end{proof}

\begin{lemma}\label{L12}
Suppose $(W,S)$ is a Coxeter system and $\Lambda$ is a graph of
groups decomposition of $W$. If $A\subset S$, $\langle A\rangle$
is virtually abelian and for any unrelated $x,y$ in $A$, $\{x,y\}$
stabilizes a vertex of $T_{\Lambda}$ (the Bass-Serre tree for
$\Lambda$), then $A$ stabilizes a vertex of $T_{\Lambda}$. In
particular, if $\Psi$ is the visual graph of groups for $\Lambda$
given by the construction for theorem \ref{MT1}, then $A\subset V$
for $\langle V\rangle$ a vertex group of $\Psi$.
\end{lemma}

\begin{proof} 
The group $\langle A\rangle$ decomposes as
$\langle a_1,b_1\rangle \times \cdots \times \langle
a_n,b_n\rangle \times \langle F\rangle$ where $A=\{a_1,b_1,\ldots
,a_n,b_n\}\cup F$, $m(a_i,b_i)=\infty$ and $F$ generates a
complete subdiagram of $\Gamma (W,S)$. If the statement of the
lemma fails, assume $n$ is minimal among all counterexamples.
Note that $n>0$ since $F$ is FA (see \cite{MT}). By the minimality
of $n$, $\langle a_1,\ldots , a_n,b_2,\ldots ,b_n,F\rangle$
stabilizes a vertex $V_1$ of $T_{\Lambda}$, $\langle a_2,\ldots ,
a_n,b_1,\ldots ,b_n,F\rangle$ stabilizes $V_2$ and $\langle
a_1,b_1\rangle$ stabilizes $V_3$. As $T_{\Lambda}$ is a tree,
there is a vertex $V$ of $T_{\Lambda}$ common to the three
geodesics connecting pairs in $\{V_1,V_2,V_3\}$, and $A$
stabilizes $V$. 
\end{proof}

\section{Weak $M(\Psi,\langle V\rangle)$-JSD Decompositions}\label{Weak}

If $(W,S)$ is a Coxeter system, $\Psi$ is a reduced visual graph
of groups decomposition for $(W,S)$ and $V\subset S$ is a vertex
of $\Psi$, then a reduced $\Psi$-compatible graph of groups
decomposition, $\Lambda$, of $\langle V\rangle$ is {\it weakly}
$M(\Psi ,\langle V\rangle)$-{\it JSJ} if
\begin{enumerate}

\item each edge group of $\Lambda$ is in $M(\Psi, \langle V\rangle)$
and

\item each element of $M(\Psi,\langle V\rangle)$ is a subgroup of a
conjugate of a vertex group of $\Lambda$.
\end{enumerate}

A reduced  $\Psi$-compatible visual graph of groups decomposition
$\Psi_1$ of $\langle V\rangle$ {\it looks weakly $M(\Psi,\langle
V\rangle )$-JSJ} if

\begin{enumerate}

\item each edge group of $\Psi_1$ is in $M(\Psi,\langle V\rangle)$
and

\item no edge group $\langle E\rangle$ ($E\subset V$) of $\Psi_1$ is
crossing.

\end{enumerate}

At this point we consider decompositions that are precursors to the $n^{th}$-stage JSJ-decompositions. For a finitely generated Coxeter system $(W,S)$ we say a visual graph of groups decomposition $\Psi$ is {\it JSJ-amenable} if $\Psi$ is reduced and for any vertex $V\subset S$ of $\Psi$ and any $E\subset S$ such that $\langle E\rangle$ is the group of an edge incident to $V$, there is no $A\subset V$ such that $\langle A\rangle$ is virtually abelian and $A$ separates $E$ in $\Gamma$. In particular, if $\Psi$ is JSJ-amenable, $V\subset S$ is a vertex of $\Psi$, $A\subset V$ separates $V$ in $\Gamma$ and $\langle A\rangle\in M(\Psi, \langle V\rangle)$ then there is no $x\in E(A)$ such that $A-\{x\}$ separates $V$ in $\Gamma$. (This remark should be compared with the one following corollary \ref{C7}). In section  \ref{JSJ} (proposition \ref{amenable}) we show that $n^{th}$-stage JSJ-decompositions are JSJ-amenable.    

\begin{proposition}\label{P13}
Suppose $(W,S)$ is a Coxeter system, $\Psi$ is a visual
graph of groups decomposition for $(W,S)$ that is JSJ-amenable and $V\subset S$ is a
vertex of $\Psi$. A reduced visual $\Psi$-compatible graph of
groups decomposition $\Psi _1$ of $\langle V\rangle$ looks weakly
$M(\Psi, \langle V\rangle)$-JSJ if and only if it is weakly
$M(\Psi, \langle V\rangle)$-JSJ.
\end{proposition}

\begin{proof} 
Assume $\Psi_1$ is weakly $M(\Psi, \langle
V\rangle )$-JSJ. If $\Psi_1$ does not look weakly $M(\Psi, \langle
V\rangle )$-JSJ, then there is $E\subset V$ such that $\langle E\rangle $ is the edge group of an edge $E'$ of $\Psi_1$
and $T\subset V$ such that $\langle T\rangle \in M(\Psi, \langle
V\rangle )$ and $\langle E\rangle$ crosses $\langle T\rangle$ in
$M(\Psi, \langle V\rangle)$. Assume $E$ separates elements $t_1$
and $t_2$ of $T$ in $\Gamma$. By lemma \ref{L10},
$E=\{e_1,e_2\}\cup N $ and $T=\{t_1,t_2\}\cup N$ where $e_1$ and
$e_2$ are unrelated, $t_1$ and $t_2$ are unrelated and $N$
commutes with $\{e_1,e_2,t_1,t_2\}$. Since $\Psi_1$ is weakly
$M(\Psi, \langle V\rangle )$-JSJ, there is a vertex $U\subset V$
of $\Psi_1$ such that $\langle T\rangle$ is conjugate to a
subgroup of $\langle U\rangle$. By proposition \ref{P5},
$\{t_1,t_2\}\subset U$. As $\Psi_1$ is a tree, we may assume $E'$
is an edge of $\Psi_1$ incident to $U$. Let $Q$ be the vertex of
$E'$ opposite $U$. Let $\Psi '$ be the graph of groups
decomposition obtained from $\Psi$ by replacing $\langle
V\rangle$ by $\Psi _1$. Since $\Psi_1$ is reduced, there exists
$x\in Q-E$. By proposition \ref{P11}, $\Gamma-E$ has exactly two
components which intersect $V$ non-trivially, one containing
$t_1$ and the other containing $t_2$. Hence $x$ can be connected
to $t_1$ or $t_2$ by a path in $\Gamma-E$. This is
impossible as $E'$ separates $U$ and $Q$ in $\Psi '$ and so by
lemma \ref{LC1}, $E$ separates $U-E$ and $Q-E$ in $\Gamma$ .

Suppose $\Psi_1$ looks weakly $M(\Psi, \langle V\rangle )$-JSJ
and $B\in M(\Psi, \langle V\rangle )$. By proposition \ref{P6},
there is $A\subset V$ such that $\langle A\rangle \in M(\Psi, \langle V\rangle )$
and some conjugate of $\langle A\rangle$ is a subgroup of (finite
index in) $B$.

Suppose $x$ and $y$ are elements of $A$ and there is no vertex
group of $\Psi_1$ containing $\{x,y\}$. Then $x$ and $y$ are
separated by $E$ in $\Gamma (W,S)$ for $E$ an edge of $\Psi_1$.
This is impossible as $\Psi_1$ looks weakly $M(\Psi, \langle
V\rangle )$-JSJ. We conclude that $\{x,y\}\subset U$ for some
$U\subset V$ a vertex of $\Psi_1$.

By lemma \ref{L12}, $A$ is a subgroup of a conjugate of a vertex
group of $\Psi_1$. Let $T$ be the Bass-Serre tree for $\Psi_1$ and
$U$ a vertex of $T$ stabilized by $\langle A\rangle$. Let $B'$ be
a conjugate of $B$ such that $\langle A\rangle$ has finite index
in $B'$. If the cosets of $\langle A\rangle$ in $B'$ are $\langle
A\rangle, b_1\langle A\rangle, \ldots ,b_n\langle A\rangle$, then
the orbit of $B'U$ in $T$ is $U, b_1U,\ldots ,b_nU$. By corollary
4.8 of \cite{DicksDunwoody}, $B'$ stabilizes some vertex of $T$.
Equivalently, $B$ is a subgroup of a conjugate of a vertex group
of $\Psi_1$, and hence $\Psi_1$ is weakly $M(\Psi, \langle
V\rangle )$-JSJ. 
\end{proof}

\begin{lemma}\label{L14}
Suppose $(W,S)$ is a Coxeter system, $\Psi$ a reduced visual
graph of groups decomposition for $(W,S)$, $V\subset S$ a
vertex of $\Psi$, $\Lambda$ a reduced $\Psi$-compatible
graph of groups decomposition of $\langle V\rangle$, $\Psi_1$
the reduced visual decomposition for $\Lambda$ from the
construction for theorem \ref{MT1}, $E\subset V$ such that $\langle E\rangle$ is an edge group of
$\Psi_1$, $T$ the Bass-Serre tree for $\Lambda$ and $E'$
the edge of $T$ such that $E=\{v\in V:v \ stabilizes \ E'\}$. If
$\{x,y\}\subset V-E$ and $x$ (respectively $y$) stabilizes the
vertex $X$ (respectively $Y$) of $T$ where $X$ and $Y$ are on
different sides of $E'$, then $x$ and $y$ are in different
components of $\Gamma -E$.
\end{lemma}

\begin{proof} 
If $\Phi$ is obtained from $\Psi$ by
replacing $V$ by $\Lambda$, then $T$ is a subtree of the
Bass-Serre tree for $\Phi$. Hence, $E$ is an edge of $\Psi_{\Phi}$
the visual decomposition for $\Phi$ given by the construction for
theorem \ref{MT1}. Now apply lemma \ref{LC1}. 
\end{proof}

\begin{proposition}\label{P15}
Suppose $(W,S)$ is a Coxeter system, $\Psi$ is a visual
graph of groups decomposition for $(W,S)$ that is JSJ-amenable, $V\subset S$ is a
vertex of $\Psi$, $\Lambda$ is a reduced graph of groups
decomposition of $\langle V\rangle$ with edge groups in $M(\Psi,
\langle V\rangle)$, and $\Psi_1$ is the reduced visual
decomposition for $\Lambda$ given by the construction for theorem
\ref{MT1}. Then $\Psi_1$ is weakly $M(\Psi, \langle V\rangle)$-JSJ
if and only if $\Lambda$ is weakly $M(\Psi, \langle
V\rangle)$-JSJ.
\end{proposition}

\begin{proof} 
By lemmas \ref{Compat} and \ref{LM11},
$\Psi_1$ and $\Lambda$ are $\Psi$-compatible, and for each $E\subset V$ such that
$\langle E\rangle$ is an edge group of $\Psi_1$, $\langle E\rangle$ is in $M(\Psi,\langle V\rangle)$
and $E$ separates $V$ in $\Gamma(W,S)$. Assume $\Lambda$ is weakly
$M(\Psi, \langle V\rangle)$-JSJ. If $\Psi_1$ is not weakly
$M(\Psi, \langle V\rangle)$-JSJ, then proposition \ref{P13}
implies $\Psi_1$ does not look weakly $M(\Psi, \langle
V\rangle)$-JSJ. I.e. there are $E$ and $B$, subsets of $V$ that
separate $\Gamma (W,S)$ and generate crossing members of $M(\Psi,
\langle V\rangle)$, such that $E$ separates (in $\Gamma$)
elements $b_1$ and $b_2$ of $B$, and $\langle E\rangle$ is an edge group of $\Psi_1$.
Let $T$ be the Bass-Serre tree for $\Lambda$. The construction of
visual decompositions for theorem \ref{MT1} implies there is an
edge $E'$ of $T$ such that $E=\{v\in V:v \ stabilizes \ E'\}$ and
since $\Psi_1$ is reduced, we may assume there are elements $x$ and $y$ of
$V-E$ that stabilize verticies $X$ and $Y$ (respectively) of $T$ on opposite sides of $E'$.
Since $\Lambda$ is weakly $M(\Psi, \langle V\rangle)$-JSJ, there
is a vertex $U$ of $T$ such that $B$ stabilizes $U$. By
proposition \ref{P11}, $\Gamma -E$ has exactly two components
which intersect $V$ non-trivially, one containing $b_1$ and the
other containing $b_2$. There is an element $z$ of $V-E$ that
stabilizes a vertex of $T$ on the side of $E'$ opposite $U$. But
by lemma \ref{L14}, $z$ cannot be in the same component of
$\Gamma -E$ as $b_1$ or $b_2$, which is nonsense. The converse is
trivial. 
\end{proof}

\section{ $M(\Psi,\langle V\rangle)$-JSJ Decompositions}
 
Suppose $(W,S)$ is a Coxeter system, $\Psi$ is a reduced visual
graph of groups decomposition for $(W,S)$, and $V\subset S$ is a
vertex of $\Psi$. We define, a weakly $M(\Psi, \langle
V\rangle)$-JSJ decomposition $\Lambda$, to be $M(\Psi, \langle
V\rangle)$-{\it JSJ} if for any vertex group $U$ of $\Lambda$ and
non-trivial $\Lambda$-compatible splitting of $U$ over an
$M(\Psi, \langle V\rangle)$ subgroup, the resulting (reduced) decomposition
of $\Lambda$ is not weakly $M(\Psi, \langle V\rangle)$-JSJ.

Suppose $\Psi_1$ is a visual weak $M(\Psi, \langle V\rangle)$-JSJ
decomposition, $E\subset V$ such that $\langle E\rangle$ is a non-crossing
element of $M(\Psi, \langle V\rangle)$, and $E$ separates $U$ in
$\Gamma (W,S)$ for $U$ a vertex of $\Psi_1$. Since $\langle
E\rangle$ is non-crossing, this splitting is $\Psi_1$-compatible,
visual and non-trivial. By lemma \ref{Compat}, the resulting
decomposition of $\Psi _1$ is compatible with $\Psi$. We say a
visual weakly $M(\Psi, \langle V\rangle)$-JSJ decomposition
$\Psi_1$ {\it looks $M(\Psi, \langle V\rangle)$-JSJ} if for any $E\subset V$ such that $\langle E\rangle$ is a non-crossing member of $M(\Psi, \langle V\rangle)$ and vertex $U$ of $\Psi_1$ such that
$E\subset U\subset V$, $E$ does not separate $U$ in $\Gamma$.

The next proposition implies the existence of visual $M(\Psi,\langle V\rangle )$-JSJ decompositions for a given JSJ-amenable graph of groups decomposition $\Psi$. 

\begin{proposition}\label{P16}
Suppose $(W,S)$ is a Coxeter system, $\Psi$ is a visual
graph of groups decomposition of $(W,S)$ that is JSJ-amenable, and $V\subset S$ is a
vertex of $\Psi$. A visual weak $M(\Psi, \langle V\rangle)$-JSJ
decomposition $\Psi_1$ of $\langle V\rangle$, looks $M(\Psi,
\langle V\rangle)$-JSJ if and only if $\Psi_1$ is $M(\Psi,
\langle V\rangle)$-JSJ.
\end{proposition}

\begin{proof}  
Suppose $\Psi_1$ looks $M(\Psi, \langle
V\rangle)$-JSJ, but is not $M(\Psi, \langle V\rangle)$-JSJ. Then
there is $D\in M(\Psi, \langle V\rangle)$ and a vertex $U\subset
V$ of $\Psi_1$ such that $\langle U\rangle$ splits non-trivially
and $\Psi_1$-compatibly as $A\ast _DB$, and the resulting reduced
decomposition $\Lambda$ (of $\langle V\rangle$), is weakly
$M(\Psi, \langle V\rangle)$-JSJ. Let $\Psi'$ be the reduced
visual decomposition for $A\ast _DB(=\langle U\rangle)$ given by
the construction for theorem \ref{MT1}. As no conjugate of
$\langle U\rangle$ is contained in $A$ or $B$, $\Psi'$ has more
than one vertex. Any visual subgroup of $\langle U\rangle$
contained in a conjugate of $A$ or $B$ is contained in a vertex
group of $\Psi'$ by the construction for theorem \ref{MT1}. Hence
$\Psi'$ is compatible with $\Psi_1$. Let $\Psi_1'$ be obtained
from $\Psi_1$ by replacing $\langle U\rangle$ by $\Psi'$, and
reducing. Then $\Psi_1'$ is a reduced visual graph of groups
decomposition for $\Lambda$ as given by the construction for
theorem \ref{MT1}. By proposition \ref{P15}, $\Psi_1 '$ is weakly
$M(\Psi, \langle V\rangle)$-JSJ and so $\Psi_1$ did not look
$M(\Psi, \langle V\rangle)$-JSJ to begin with. The converse is
trivial. 
\end{proof}

\begin{theorem}\label{T17}
Suppose $(W,S)$ is a Coxeter system, $\Psi$ is a
visual graph of groups decomposition of $(W,S)$ that is JSJ-amenable, and $V\subset S$
is a vertex of $\Psi$. If $\Lambda$ is a reduced $M(\Psi, \langle
V\rangle)$-JSJ graph of groups decomposition of $\langle
V\rangle$, and $\Psi_1$ is the reduced visual decomposition
derived from $\Lambda$ by the construction for theorem \ref{MT1},
then
\begin{enumerate}
\item  $\Psi_1$ is an $M(\Psi, \langle V\rangle)$-JSJ decomposition.

\item  There is a (unique) bijection $\alpha$ of the vertices
of $\Lambda$ to the vertices of $\Psi_1$ such that for each vertex
$U$ of $\Lambda$, $\Lambda(U)$ is conjugate to $\Psi(\alpha (U))$.

\item  Each edge group of $\Lambda$ is conjugate to a special
subgroup of $\langle V\rangle$.
\end{enumerate}
\end{theorem}

\begin{proof}
The decomposition $\Psi_1$ is weakly
$M(\Psi, \langle V\rangle)$-JSJ by Proposition \ref{P15}. Suppose
$U$ is a vertex of $\Lambda$ with vertex group $A=\Lambda (U)$. We
wish to show that $A$ is a subgroup of a conjugate of a vertex
group of $\Psi_1$. Otherwise, the action of $A$ on $T$, the
Bass-Serre tree for $\Psi_1$, defines a non-trivial reduced graph
of groups decomposition $\Phi$ of $A$ such that each edge group
of $\Phi$ is a subgroup of a conjugate of an ($M(\Psi, \langle
V\rangle)$) edge group of $\Psi_1$. If $A$ does not stabilize a
vertex of $T$, then $\Phi$ is nontrivial.  If $C$ is an edge
group of $\Lambda$ incident with $A$, then $C$ contains a
conjugate of an edge group $Q$ of $\Psi_1$. Since $C$ and $Q$ are
in $M(\Psi, \langle V\rangle)$, a conjugate of $Q$ is of finite
index in $C$ and at the same time stabilizes an edge of $T$.  But
then the orbit of this edge in $T$ under the action of $C$ is
finite.  By Corollary 4.8 of \cite{DicksDunwoody}, $C$ stabilizes
some vertex of $T$ and so is contained in a conjugate of a vertex
group of $\Phi$. Hence $\Phi$ is compatible with $\Lambda$. We
also require $\Phi$ to be compatible with $\Psi$. Suppose
$D\subset V$ is an edge of $\Psi$ incident to $V$, such that some
conjugate of $\langle D\rangle$ is a subgroup of $A$. Then
$D\subset Q$ for $Q$ a vertex of $\Psi _1$. Hence $\langle
D\rangle$ stabilizes a vertex group of $T$ and so is a subgroup
of a conjugate of a vertex group of $\Phi$, as required.
Replacing $A$ in $\Lambda$ by this decomposition and reducing
gives a $\Psi$-compatible reduced graph of groups decomposition
$\Lambda '$ of $\langle V\rangle$ with $M(\Psi, \langle
V\rangle)$ edge groups. By hypothesis, $\Lambda '$ is not weakly
$M(\Psi, \langle V\rangle)$-JSJ. Hence there exists $B$, an
$M(\Psi, \langle V\rangle)$ subgroup of $\langle V\rangle$ such
that $\langle V\rangle$ splits non-trivially over $B$ and $B$ is
not a subgroup of a conjugate of a vertex group of $\Lambda '$.
Since $\Lambda$ is weakly $M(\Psi, \langle V\rangle)$-JSJ, $B$ is
a subgroup of a conjugate of $A$. But, since $\Psi_1$ is weakly
$M(\Psi, \langle V\rangle)$-JSJ, $B$ is a subgroup of a conjugate
of a vertex group of $\Psi_1$ and hence of a vertex group of
$\Phi$ and of $\Lambda '$. We conclude that $A$ is a subgroup of
a conjugate of a vertex group of $\Psi_1$, as desired.

If $A(=\Lambda (U))$ is a subgroup of a conjugate of
$\Psi_1(U')(\equiv \langle U'\rangle)$ for $U'$ a vertex of
$\Psi_1$, then, since $\Psi_1(U')$ is a subgroup of a conjugate
of a vertex group of $\Lambda$, $A$ is a subgroup of a conjugate
of a vertex group of $\Lambda$. By lemma \ref{Conj} the vertex
group $A$ at $U$ is a subgroup of a conjugate of a vertex group
at $U''$ only if $U=U''$ and (since $\Lambda$ is a tree) the
conjugate is by an element of $A$. But then $A$ is conjugate to
$\Psi_1(U')$. Since no vertex group of $\Psi_1$ is contained in a
conjugate of another, $U'$ is uniquely determined, and we set
$\alpha(U)=U'$. Since each vertex group $\Psi_1(U')$ is contained
in a conjugate of some $\Lambda(U)$ which is in turn conjugate to
$\Psi_1(\alpha(U))$ we must have $U'=\alpha(U)$ and each $U'$ is
in the image of $\alpha$.

If $\Psi_1$ is not $M(\Psi, \langle V\rangle)$-JSJ, then it does
not look $M(\Psi, \langle V\rangle)$-JSJ and some vertex group
$W_1$ of $\Psi_1$ visually splits nontrivially and $\Psi _1$-
compatibly over an $M(\Psi, \langle V\rangle)$ visual subgroup
$U_1$ to give a $\Psi$-compatible visual graph of groups
decomposition $\Phi$ of $\langle V\rangle$ with $U_1$ an edge
group. Now $W_1$ is conjugate to a vertex group $A$ of
$\Lambda$.  As a subgroup of $\langle V\rangle$, $A$ acts on the
Bass-Serre tree $T'$ for $\Phi$, but $A$ cannot stabilize a
vertex of $T'$, otherwise $W_1$ stabilizes a vertex of $T'$ and
we assumed $W_1$ split nontrivially. As above this contradicts
the assumption that $\Lambda$ is $M(\Psi, \langle V\rangle)$-JSJ,
implying instead that $\Psi_1$ is $M(\Psi, \langle V\rangle)$-JSJ.

Since $\Lambda$ is a tree, we can take each edge group of
$\Lambda$ as contained in its endpoint vertex groups taken as
subgroups of $\langle V\rangle$. Hence each edge group is simply
the intersection of its incident vertex groups (up to
conjugation). Since vertex groups of $\Lambda$ correspond to
conjugates of vertex groups in $\Psi_1$, their intersection is
conjugate to a visual subgroup by Lemma \ref{L3}. 
\end{proof}

Theorem \ref{T17} shows that all $M(\Psi, \langle V\rangle)$-JSJ decompositions for Coxeter
groups are basically visual. The next collection of lemmas lead to a proof of the uniqueness
of $M(\Psi, \langle V\rangle)$-JSJ decompositions for Coxeter
groups.

\begin{lemma}\label{L18}
Suppose $(W,S)$ is a Coxeter system, $\Psi$ is a visual
graph of groups decomposition of $(W,S)$ that is JSJ-amenable, $V\subset S$ is a
vertex of $\Psi$ and $\Psi_1 $ is a visual $M(\Psi,\langle
V\rangle)$-JSJ decomposition for $\langle V\rangle$. If $D\subset
V$ is such that $D$ separates $V$ in $\Gamma (W,S)$ and $\langle
D\rangle\in M(\Psi,\langle V\rangle)$, then there is $Q\subset V$
a vertex of $\Psi_1$ such that $D\subset Q$. Furthermore, if
$\langle D\rangle$ is non-crossing in $M(\Psi, \langle V\rangle)$, then $E(D)$ is a subset of an
edge group of $\Psi_1$.
\end{lemma}

\begin{proof} 
By definition some conjugate of $D$ is a
subset of $\langle Q\rangle $ for $Q$ a vertex of $\Psi_1$, but
we require more. By proposition \ref{P5}, $E(D)\subset
Q$. We have $D=E(D)\cup N$ where $N\subset V$ generates a finite
subgroup of $\langle V\rangle$ which commutes with $E(D)$. As $N$
determines a complete subdiagram of $\Gamma (\langle V\rangle,V)$, $N$ is a
subset of some vertex of $\Psi_1$. Let $D'$ be a largest subset
of $D$ such that $N\subset D'$, $D'\subset Q_1\subset V$ and
$Q_1$ is a vertex of $\Psi_1$. If $D'\ne D$, let $a\in D-D'$.
Observe that $\{a\}\cup N$ generates a finite group and so there
exists $Q_2(\subset V)$ a vertex of $\Psi_1$ such that $\{a\}\cup
N\subset Q_2$. If $U$ is a vertex common to the three
$\Psi_1$-geodesics connecting pairs in $\{Q,Q_1,Q_2\}$, then
$\{a\}\cup D'\subset U$ contrary to the definition of $D'$. We
conclude $D'=D$ as required.
 
The second part of the lemma is not used in the rest of the paper and is tedious to prove. The reader may wish to skip this argument on a first reading.

Assume $\langle D\rangle$ is non-crossing in $M(\Psi,\langle
V\rangle)$ and $Q$ is a vertex of $\Psi_1$ with $D\subset Q\subset
V$. Assume $E(D)\ne\emptyset$. If $Q'$ is a vertex of $\Psi_1$,
distinct from $Q$ and $E(D)\subset Q'$ then for any $F\subset V$ such that $\langle F\rangle$ is the edge group of an edge on
the $\Psi_1$-geodesic connecting $Q$ and $Q'$, $E(D)\subset F$.
Hence we may assume:

\medskip

\noindent $(\ast )$ $Q$ is the only vertex of $\Psi_1$ such that
$E(D)\subset Q$.

\medskip

Even as $D$ separates $V$ in $\Gamma$, $D$ does
not separate $Q$ in $\Gamma$ (otherwise, since $D$ is non-crossing
and does not separate $E$ (for $\langle E\rangle$ an edge group of  $\Psi$) in $\Gamma$,
$\langle Q\rangle$ splits non-trivially and compatibly with
$\Psi_1$ and $\Psi$ over $D$ - contrary to our JSJ assumption on
$\Psi_1$). Hence $Q-D\subset C$ for $C$ a component of $\Gamma
-D$. 

Next we show that if $C'$ is  a component of $\Gamma -D$ other than $C$ such that  $C'$
intersects $V$ non-trivially, and $x\in E(D)$ then there is a vertex $V_{x,C'}$ of $\Psi_1$ such that $x\in V_{x,C'}$ and $V_{x,C'}\cap C'\ne \emptyset$. First observe that $C'\cap Q=\emptyset$.  By lemma \ref{L8}, for each $x\in
E(D)$ there is an edge $[xa]$ of $\Gamma$ such that $a\in C'$. Let $x\equiv x_0,x_1,\ldots ,x_n$ be the consecutive vertices of a path in $\Gamma $ from $x$ to $x_n\in V\cap C'$ such that for $i>0$, $x_i\in C'$. 
We may assume $n$ is the smallest integer such that $x_n\in V\cap C'$. Since $x\not\in C'$, $n\ne 0$. Let $\Psi'$ be the visual graph of groups decomposition of $W$ obtained by replacing the vertex $V$ of $\Psi$ by $\Psi_1$. By lemma \ref{visgog}, for each $i\in \{1,2,\ldots ,n\}$ there is a vertex $V_i\subset S$ of $\Psi'$ such that $\{x_{i-1},x_i\}\subset V_i$. 
Define $V_0\equiv Q$. Since $Q\cap C'=\emptyset$, $V_i\ne Q$ for $i>0$. Let $\alpha _i$ be the $\Psi'$-geodesic from $Q\equiv V_0$ to $V_i$ and $\beta _i$ be the $\Psi'$-geodesic from $V_i$ to $V_{i+1}$ for $i\geq 1$. 
By lemma \ref{visgog}, if $B\subset S$ is a vertex of $\beta _i$, then $x_i\in C'\cap B$ and so $\beta_i$ does not pass through $Q$. This implies that $\alpha _n$ can be written as a non-trivial subpath $\tau$ of $\alpha _1$ followed by a path $\lambda$ where each edge of $\lambda$ is an edge of some $\beta _i$. 
Let $Q\equiv X_1,\ldots X_k\equiv V_n$ be the consecutive vertices of $\alpha _n$. If $X_i$ is the end point of $\tau$ then $x_0\in X_i$ (lemma \ref{visgog} implies $x_0$ is an element of every vertex of $\alpha$). Since $X_i$ is also the initial point of $\lambda$, $\{x_0,x_m\}\subset X_i$ for some $m\in \{1,\ldots ,n\}$. Since $\Psi_1$ is a subtree of $\Psi'$, there is $j\in \{1,\ldots ,k\}$ such that $X_1,\ldots ,X_j$ are vertices $\Psi_1$ and $X_{j+1},\ldots ,X_k$ are not. 
If $j\geq i$, then as $\{x_0,x_m\}\subset X_i$, set $V_{x,C'}=X_i$ to finish. If $j<i$, then as $\Psi_1$ is a subtree of $\Psi'$, any $\Psi'$-geodesic from $V_n\equiv X_k$ to a vertex of $\Psi_1$ must have initial segment with consecutive vertices $V_n\equiv X_k, X_{k-1},\ldots ,X_j$. 
Hence if $Q'\subset V$ is a vertex of $\Psi_1$ containing $x_n(\in V\cap C')$, then the $\Psi'$ geodesic from $V_n$ to $Q'$ passes through $X_j$. As $x_n\in V_n\cap Q'$, lemma \ref{visgog} implies $x_n\in X_j$. Since $j<i$, $x_0\in X_j$ and set $V_{x,C'}=X_j$ to finish the claim.

Note that for all $x\in E(D)$ and for all components $C'\ne C$ of $\Gamma-D$, $V_{x,C'}\ne Q$ (as $C'\cap Q=\emptyset$). 
For $x\in E(D)$ and $C'\ne C$ a component of $\Gamma -D$, there may be more than one possible choice for $V_{x,C'}$. Assume that $E(D)\cap V_{x,C'}$ is maximal over all possible choices. I.e. if $X$ is a vertex of $\Psi_1$, $x\in X$, and $X\cap C'\ne\emptyset$ then the number of elements of $X\cap E(D)$ is less that or equal to the number of elements of $V_{x,C'}\cap E(D)$. 
By ($\ast$), $E(D)\not \subset V_{x,C'}$.
Let $y\in E(D)-V_x$, then $Q$, $V_{y,C'}$ and $V_{x,C'}$ are distinct.
There is a vertex $U$ common to the three geodesics of $\Psi_1$
connecting pairs in $\{Q,V_{x,C'},V_{y,C'}\}$ and so $[E(D)\cap
V_{x,C'}]\cup\{y\}\subset U$. In particular, $x\in U$. By the maximality condition on $V_{x,C'}\cap E(D)$, $U\cap
C'=\emptyset$. Hence $V_{x,C'}\ne U\ne V_{y,C'}$ and $U$ separates $V_{x,C'}$
and $V_{y,C'}$ in $\Psi_1$. If $\Psi'$ is the graph of groups decomposition of $W$ obtained from $\Psi$ by replacing $V$ by $\Psi_1$, then $\Psi_1$ is a subtree of the tree $\Psi'$ and $U$ separates $ V_{x,C'}$ and $V_{y,C'}$ in $\Psi'$. Suppose $\alpha$ is an edge path in $C'$
(with consecutive vertices $c_0,\ldots ,c_n$) connecting $c_0\in
C'\cap V_{x,C'}$ to $c_n\in C'\cap V_{y,C'}$. Let $C_0=V_{x,C'}$, $C_i$ be a
vertex of $\Psi '$ such that $C_i$ contains $\{c_{i-1},c_i\}$ and
$C_{n+1}=V_{y,C'}$. If $X$ is a vertex of the $\Psi'$-geodesic from
$C_i$ to $C_{i+1}$, then lemma \ref{visgog} implies $X$ contains $c_i\in C'$. Hence, $\alpha$
defines a path in $\Psi'$ from $V_{x,C'}$ to $V_{y,C'}$ avoiding $U$
(which is impossible). 
\end{proof}

For a Coxeter system $(W,S)$ and vertex $V\subset S$
of a visual $M(W,S)$-JSJ decomposition of $(W,S)$, $V$ may not be
a connected subset of $\Gamma (W,S)$.

\begin{figure}[ht]
\begin{center}
\epsfig{file=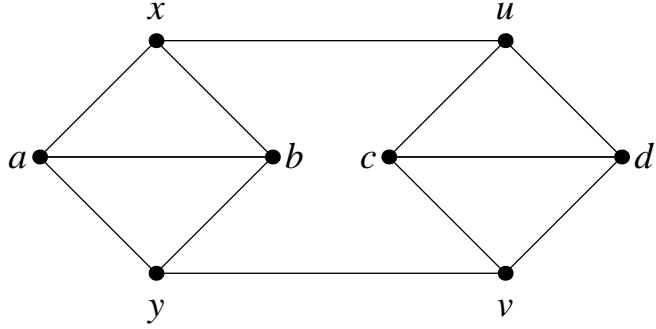}
\caption{A disconnected JSJ vertex group \label{jsjfig3}}
\end{center}
\end{figure}

\begin{example}
Let $(W,S)$ be the Coxeter system given by $\Gamma (W,S)$ pictured in figure \ref{jsjfig3}, where
each edge has label 3. The only non-crossing visual virtually abelian splitting subgroups for
this system are $\langle x,y\rangle$ and $\langle u,v\rangle$ ($\langle x,v\rangle$ and $\langle u,y\rangle$ are crossing).
The JSJ-decomposition is given by: $$ \langle a,b,x,y\rangle \ast
_{\langle x,y\rangle} \langle x,y,u,v\rangle \ast _{\langle
u,v\rangle} \langle u,v,c,d\rangle $$ The set $\{x,y,u,v\}$
generates a vertex group of this decomposition and is not
connected in $\Gamma$. 
\end{example}

Still we have the following:

\begin{lemma}\label{L19}
Suppose $(W,S)$ is a Coxeter system, $\Psi$ is a visual
graph of groups decomposition of $(W,S)$ that is JSJ-amenable, $V\subset S$ is a
vertex of $\Psi$, and $\Psi_1$ is a visual $M(\Psi ,\langle
V\rangle)$-JSJ decomposition of $\langle V\rangle$.
If $D\subset V$, $\langle D\rangle\in M(\Psi, \langle V\rangle)$,
$D$ separates $V$ in $\Gamma (W,S)$, and $\langle D\rangle$ is
non-crossing in $M(\Psi,\langle V\rangle)$, then there is no
$Q\subset V$ a vertex of $\Psi_1$ such that $D$ separates $Q$ in
$\Gamma(W,S)$.
\end{lemma}

Note that in the previous example, neither $\{x,y\}$ nor
$\{u,v\}$ separates $\{x,y,u,v\}$ in $\Gamma$.

\medskip

\begin{proof} 
Suppose otherwise. By lemma \ref{L18},
there is $Q'\subset V$ a vertex of $\Psi_1$ such that $D\subset
Q'$. As $\Psi_1$ is $M(\Psi,\langle V\rangle)$-JSJ, $D$ does not separate $Q'$ in
$\Gamma$. In particular, $Q\ne Q'$. Suppose $F\subset V$ is such that $\langle F\rangle$ is the group of an
edge of the $\Psi_1$-geodesic connecting the vertices $Q$ and
$Q'$. As $F$ separates $Q$ and $Q'$ in the tree $\Psi_1$, (and
since $D(\subset Q'$) separates $Q$ in $\Gamma$) $F$ separates
$Q$ in $\Gamma$. But if $\langle F\rangle$ is the group of the edge of this geodesic incident
to $Q$, there is an induced splitting of $\langle Q\rangle$ over
$\langle F\rangle$ compatible with $\Psi_1$ (since $\langle F\rangle$ is non-crossing) and $\Psi$ (since $\Psi$ is JSJ-amenable). This is
impossible as $\Psi_1$ is $M(\Psi, \langle V\rangle)$-JSJ.
\end{proof}

\begin{theorem}\label{T20}
Suppose $(W,S)$ is a Coxeter system, $\Psi$ is a visual
graph of groups decomposition of $(W,S)$ that is JSJ-amenable, $V\subset S$ is a
vertex of $\Psi$, and $\Psi_1$ and $\Psi_2$ are visual $M(\Psi
,\langle V\rangle)$-JSJ decomposition of $\langle V\rangle$. If $Q\subset V$ is a vertex of
$\Psi_1$, then $Q$ is a vertex of $\Psi_2$. I.e. The vertex groups of two visual
$M(\Psi, \langle V\rangle)$-JSJ decompositions of $\langle
V\rangle$ are the same.
\end{theorem}

\begin{proof} 
It suffices to show that $Q\subset Q'$, for
$Q'$ a vertex of $\Psi _2$. Otherwise, there
exists $D\subset V$ an edge of $\Psi_2$ and $D$ separates $Q$ in
$\Gamma$. This is contrary to lemma \ref{L19}. 
\end{proof}

Theorems \ref{T17} and \ref{T20} imply:

\begin{corollary} \label{C21}
Suppose $(W,S)$ is a Coxeter system, $\Psi$ is a visual
graph of groups decomposition of $(W,S)$ that is JSJ-amenable, $V\subset S$ a vertex
of $\Psi$, $\Lambda_1$ and $\Lambda_2$ are $M(\Psi ,\langle
V\rangle)$-JSJ decomposition of $\langle V\rangle$. Then there is a (unique) bijection $\alpha$ of
the vertices of $\Lambda _1$ to the vertices of $\Lambda _2$ such
that for each vertex $Q$ of $\Lambda _1$, $\Lambda _1(Q)$ is
conjugate to $\Lambda _2(\alpha (Q))$. $\square$
\end{corollary}

\begin{remark} For a Coxeter system $(W,S)$, the vertex groups of
two visual $M(W,S)$-JSJ decompositions must be the same, but it is
unreasonable to expect two such graph of groups decompositions to
be exactly the same. As an example, a Coxeter diagram $\Gamma
(W,S)$ may be such that $S=A\cup B\cup C$ where $A\cap B=A\cap
C=B\cap C=E$ and a visual $M(W,S)$-JSJ decomposition is given by
$\langle A\rangle \ast _{\langle E\rangle }\langle B\rangle \ast
_{\langle E \rangle }\langle C\rangle$ or $\langle A\rangle \ast
_{\langle E\rangle }\langle C\rangle \ast _{\langle E \rangle
}\langle B\rangle$. See for instance example \ref{E4} of section
\ref{Orbi}.
\end{remark}

\section{JSJ-Decompositions}\label{JSJ}

Suppose $(W,S)$ is a finitely generated Coxeter system. Let
$\Psi_0$ be the trivial (single vertex) graph of groups
decomposition of $W$. Then $\Psi_0$ is vacuously JSJ-amenable. If $\Psi_1$ is the (unique) visual
$M(W,S)(=M(\Psi_0,\langle S\rangle))$-JSJ decomposition of $W$ and $E\subset S$ is such that $\langle E\rangle$ is the group of an edge incident to $V\subset S$ for $V$ a vertex of $\Psi_1$, then $\langle E\rangle\in M(\Psi_0,\langle S\rangle)$. 
If $A\subset V$ is such that $\langle A\rangle$ is virtually abelian and $A$ separates $E$ in $\Gamma(W,S)$ then lemma \ref{L10} implies $\langle A\rangle\in M(\Psi_0,\langle S\rangle)$. 
But then $\langle E\rangle$ and $\langle A\rangle$ are crossing in $M(\Psi_0,\langle S\rangle)$, contrary to the fact that $\langle E\rangle$ is an edge group of $\Psi_1$ (the (unique) visual $M(W,S)(=M(\Psi_0,\langle S\rangle))$-JSJ decomposition of $W$). Instead, $\Psi_1$ is JSJ-amenable.

If $\langle V\rangle$ is a vertex group of $\Psi_1$, lemma
\ref{Compat} implies the (visual) $M(\Psi_1,\langle V\rangle)$-JSJ
decomposition of $\langle V\rangle$ is compatible with $\Psi_1$.
Inductively assume $\Psi_j$ is JSJ-amenable for all $j<i$. Let $\Psi_i$ be obtained from $\Psi _{i-1}$ by
replacing each vertex group $\langle V\rangle$ of $\Psi_{i-1}$ by
the (unique) $M(\Psi_{i-1},\langle V\rangle)$-JSJ decomposition
insured by theorem \ref{T20}. By corollary \ref{C21}, we may
assume this decomposition is visual. The compatibility of this
decomposition with $\Psi_{i-1}$ is insured by repeated application
of lemma \ref{Compat}. Our first result in this section will be to show $\Psi_i$ is JSJ-amenable for all $i$, insuring our definitions are meaningful. We call $\Psi _i$ the {\it $i^{th}$-stage
of the JSJ-decomposition for $(W,S)$}. As all decompositions
involved are visual, there is an integer $n$ such that for every
vertex group $\langle V\rangle$ of $\Psi_n$, and $\langle
D\rangle \in M(\Psi_n,\langle V\rangle$) if $D$ separates
$\langle V\rangle$ in $\Gamma(W,S)$ then $\langle D\rangle$ is
crossing. In this case $\Psi_n=\Psi_{n+1}$ and we define $\Psi_n$
to be a {\it JSJ-decomposition of $W$}. This decomposition is
unique in the sense of corollary \ref{C21}. In particular, if
$(W,S)$ and $(W,S')$ are finitely generated Coxeter systems, then
the JSJ-decompositions of $W$ derived from $(W,S)$ and $(W,S')$
have conjugate vertex groups.

\begin{proposition}\label{amenable}
Suppose $(W,S)$ is a finitely generated Coxeter system. Then $\Psi_n$, the $n^{th}$-stage of the JSJ-decomposition for $(W,S)$ is JSJ-amenable for all $n$.
\end{proposition}

\begin{proof}
We have shown $\Psi_0$ and $\Psi_1$ are JSJ-amenable. Assume $\Psi_{j}$ is JSJ-amenable for all $j<n$. Then $\Psi_n$ is a visual graph of groups and each edge group of $\Psi_n$ is a member of $M(\Psi_j, \langle Q\rangle)$ for $j<n$. Suppose $V\subset S$ is a vertex of $\Psi_n$ and $E\subset S$ is such that $\langle E\rangle$ is the group of an edge of $\Psi_n$ incident to $V$. We must show there is no $A\subset V$ such that $\langle A\rangle$ is virtually abelian and $A$ separates $E$ in $\Gamma$. Suppose otherwise. We may assume $\langle E\rangle$ is the group of an edge of the $M(\Psi_j,\langle U\rangle)$-JSJ decomposition of $\langle U\rangle$ for $U$ a vertex of $\Psi_j$, where $j<n$ and $V\subset U$. 
If $j<n-1$, then $A$ cannot separate $E$ in $\Gamma$ since $\Psi_j$ is JSJ-amenable. If 
$ j=n-1$, and $A$ separates $E$ in $\Gamma$, 
then lemma \ref{L10} implies $A\in M(\Psi_{n-1}, \langle U\rangle)$ and that $\langle A\rangle$ and $\langle E\rangle$ are crossing in $M(\Psi_{n-1},\langle U\rangle)$. 
But then $\langle E\rangle$ is not the group of the edge of the $M(\Psi_{n-1},\langle U\rangle)$ decomposition of $\langle U\rangle)$. Instead, $\Psi_n$ is JSJ-amenable.
\end{proof}

\begin{proposition}\label{nosplit}
Suppose $(W,S)$ is a Coxeter system, $\Psi_i$ is the $i^{th}$-stage of the JSJ-decomposition for $(W,S)$, $V\subset S$ is a
vertex of $\Psi_i$, $\langle A\rangle\in {\cal C}(W,S)$ and
$B\subset V$ such that $B$ separates $V$ in $\Gamma (W,S)$ and
$B\in M(\Psi_i, \langle V \rangle)$. If $A$ separates $B$ in
$\Gamma (W,S)$ then $A\subset V$.
\end{proposition}

\begin{proof}
This result is trivial for $i=0$ as $V=S$
in this case. Suppose the proposition fails. Let $i\geq 1$ be the
largest integer such that for all $j<i$, there is no $A\subset S$
and vertex $V$ of $\Psi_{j}$ such that $\langle A\rangle\in {\cal
C}(W,S)$, $A$ separates $B$ in $\Gamma$ for some $\langle B\rangle\in
M(\Psi_{j},\langle V\rangle)$, $B$ separates $V$ in $\Gamma$, and
$A\not \subset V$. Then there is $\langle A\rangle \in {\cal
C}(W,S)$, $V'$ a vertex of $\Psi_i$, and $B\subset V'$ such that
$B$ separates $V'$ in $\Gamma$, $\langle B\rangle \in
M(\Psi_i,\langle V'\rangle )$, $A$ separates $B$ in $\Gamma$, and
$A\not \subset V'$. Assume $\{b_1,b_2\}\subset B$ and $A$
separates $b_1$ from $b_2$ in $\Gamma$.  By lemma \ref{L10},
$B-\{b_1,b_1\}=A-\{a_1,a_2\}\equiv M$ where $B$ separates $a_1$
and $a_2$ in $\Gamma$.  Let $V$ be the vertex group of
$\Psi_{i-1}$ containing $V'$, and $\Psi _{i}'$ the visual $M(\Psi
_{i-1},\langle V\rangle)$-JSJ decomposition that splits $\langle V\rangle$ to
give $\Psi_i$.

We show $B\not \in M(\Psi _{i-1},\langle V\rangle)$. Assume otherwise. As $A$
separates $B$ in $\Gamma$ the definition of $i$ implies that
$A\subset V$. Lemma \ref{L10} implies $A\in M(\Psi_{i-1},\langle
V\rangle)$. As $M\cup \{a_1,a_2\}=A\not\subset V'$ and $M\subset B\subset V'$, we assume $a_1\not \in V'$.
Then as $\Psi _i'$ is a tree, there is $D\subset V$ such that $\langle D\rangle$ is the group of an edge of $\Psi _i'$ incident to $V'$ such
that any path in $\Gamma $ from $a_1$ to $V'$ intersects $D$.
Thus, $M\subset D$. Any path from $a_1$ to $a_2$ intersects
$B\subset V'$, and so such a path intersects $D$. As $D\subset
V'$, $a_1\not\in D$. If $a_2\not\in D$, then $D$ separates $a_1$
and $a_2$ in $\Gamma$, so that $D$ and $A$ are crossing in
$M(\Psi _{i-1}, \langle V\rangle)$. But no edge of $\Psi_i'$  is
crossing in $M(\Psi_{i-1},\langle V\rangle )$. Instead, $a_2\in
D\subset V'$. As $\{b_1,b_2\} \subset V'$, and since any path
from $b_1$ to $b_2$ contains $a_1$ or $a_2$, any path from $b_1$
to $b_2$ intersects $D$. Hence, as $\langle D\rangle$ is
non-crossing in $M(\Psi _{i-1},\langle V\rangle)$, and since $B\in
M(\Psi_{i-1},\langle V\rangle)$, $\{b_1,b_2\}\cap D\ne\emptyset$.
Furthermore, $\{b_1,b_2\}\not \subset D$, since otherwise, $A$
would cross (the non-crossing) $D$ in $M(\Psi _{i-1},\langle V\rangle)$. We assume
$b_1\in D$ and $b_2\in V'-D$. The set $\{a_2\}\cup M$ does not
separate $V$ in $\Gamma$ as $A\in M(\Psi _{i-1},\langle
V\rangle)$. Let $\alpha $ be a shortest path in $\Gamma$ from
$b_1$ to $b_2$, avoiding $\{a_2\}\cup M$. Then $\alpha$ contains
$a_1$, and so the segment of $\alpha$ from $a_1$ to $b_2$
contains $d\in D$. The segment of $\alpha$ from $d$ to $b_2$
avoids $A$. Now $\{b_1,d\}\subset D-(\{a_2\}\cup M)$. There is no
edge connecting $b_1$ and $d$, since $b_1$ cannot be connected to
$b_2$ avoiding $A$. As $\langle D\rangle$ is virtually abelian, theorem \ref{T1} implies $a_2$ is connected to
$b_1$ and to $d$ by an edge (labeled 2). The segment of $\alpha$
from $a_1$ to $d$ followed by the edge from $d$ to $a_2$ avoids
$B$, which is nonsense. Instead, $\langle B\rangle \not \in M(\Psi
_{i-1},\langle V\rangle)$ as desired.

Proposition \ref{P6} implies there is $K\subset V$ such that $K$ separates $V$ in $\Gamma$ and $\langle K\rangle$ is a maximal rank element of $M(\Psi _{i-1},\langle V\rangle )$ such that
$E(K)\subset E(B)$. Then $\langle E(K)\rangle$ has infinite
index in $B$. As above we assume $a_1\not\in V'$ and there is $D\subset V$ such that $\langle D\rangle$ is the group of an edge of $\Psi _i'$ incident to $V'$ such that any path in
$\Gamma $ from $a_1$ to $V'$ intersects $D$. In particular
$M\subset D$.

We show  $E(K)=E(M)=E(D)$. If $\{b_1,b_2\}\subset K$, then $A$
separates $K$ and lemma \ref{L10} implies
$Rank(K)=Rank(A)(=Rank(B))$ which is nonsense. Instead,
$E(K)=E(K\cap B)\subset E(M)$. As $M\subset D$, $E(K)\subset E(M)\subset E(D)$. Since it is also true that
$\{K,D\}\subset M(\Psi _{i-1},\langle V\rangle)$, we see that $E(D)= E(K)$
and so $E(K)=E(M)=E(D)$.

Next we show $E(M)$ is not a subset of the group of an edge of $\Psi _{i-1}$
(incident to $V)$. Otherwise, there is a $j<i-1$, a vertex $\hat
V$ of $\Psi_j$ containing $V$, and a non-crossing $F\in M(\Psi
_j, \langle \hat V\rangle)$ such that $F\subset V$, and
$E(M)\subset F$. But then $K$ is non-crossing in $M(\Psi
_j,\langle \hat V\rangle)$ and as $K$ separates $V$ in $\Gamma$,
the visual $M(\Psi _j\langle \hat V\rangle )$-JSJ decomposition
of $\langle \hat V\rangle$ was not maximal - which is nonsense.

As $E(M)$ is not a subset of the group of an edge of $\Psi _{i-1}$ incident to
$V$, $A\subset V$. If $a_2\not \in V'$, then there exists $D'\subset V$ such that$\langle D'\rangle$ is the group of 
an edge of $\Psi _i'$, $D'$ separates $a_2$ from $V'$ in $\Gamma$, and
$E(M)=E(D')$. Note that $D-E(M)$ and $D'-E(M)$ generate finite
groups. In particular, these sets define complete subgraphs of
$\Gamma$. As $\langle B\rangle \in M(\Psi _i, V')$, $E(M)$ does
not separate $b_1$ from $b_2$ in $\Gamma$. If $\alpha $ is a
shortest path in $\Gamma$ from $b_1$ to $b_2$ avoiding $E(M)$,
then $a_1$ or $a_2$ is a vertex of $\alpha$. If $a_1$ is a vertex
of $\alpha$, then some vertex $t$ of $\alpha$, between $b_1$ and
$a_1$ belongs to $D-E(M)$ and some vertex $s$ of $\alpha$ between
$a_1$ and $b_2$ belongs to $D-E(M)$. But there is an edge between
$t$ and $s$ contradicting the minimality of $\alpha$. Similarly if
$a_2$ is a vertex of $\alpha$. Instead, either $a_1$ or $a_2$ is
an element of $V'$. If $a_2\in V'$, then since $B\in M(\Psi
_i,\langle V'\rangle)$, $\{a_2\}\cup M$ does not separate $V'$ in
$\Gamma$. Let $\alpha$ be a shortest path from $b_1$ to $b_2$ in
$\Gamma$, avoiding $\{a_2\}\cup M$. Then $a_1$ is a vertex of
$\alpha$  and as before, there exists $t$ and $s$ in $D-E(M)$
such that $t$ is between $b_1$ and $a_1$ and $s$ is between $a_1$
and $b_2$. The edge joining $t$ and $s$ shortens $\alpha$ which
is impossible. 
\end{proof}

\begin{theorem} \label{T33}
Suppose $\Psi$ is the visual JSJ-decomposition for $(W,S)$
$\langle V\rangle$ is a vertex group of $\Psi$, $\Phi$ is an
arbitrary graph of groups decomposition of $W$ with virtually
abelian edge groups, and $T$ is the Bass-Serre tree for $\Phi$. If
$\langle A\rangle \in {\cal C}(W,S)$, then
\begin{enumerate}

\item $A$  does not separate $B$ in $\Gamma (W,S)$ if $\langle
B\rangle$ is an edge group of $\Psi$,
\item the decomposition of $\langle V\rangle$ induced by its action
on $T$ is compatible with $\Psi$, and
\item if $M(\Psi,\langle V\rangle)$ contains no crossing elements,
then $\langle V\rangle$ stabilizes a vertex of $\Phi$.
\end{enumerate}
\end{theorem}

\begin{proof}
If 1) fails, there is an integer $i$ and a
vertex group $\langle V\rangle$ of $\Psi_i$, the $i^{th}$-level
visual JSJ-decomposition for $(W,S)$, such that $\langle
B\rangle$ is an edge group of $\Psi_{i+1}'$ the $M(\Psi
_i,\langle V\rangle)$ visual JSJ-decomposition of $\langle
V\rangle$. By proposition \ref{nosplit}, $A\subset V$ and so by
lemma \ref{L10}, $\langle A\rangle$ and $\langle B\rangle$ are
crossing in $M(\Psi_i,\langle V\rangle)$. But this is impossible
as no edge of $\Psi_{i+1}'$ is crossing in $M(\Psi_i,\langle
V\rangle)$.

For the second part of the theorem, assume $B$ is an edge of
$\Psi$ and $\langle B\rangle$ does not stabilize a vertex of $T$
(equivalently, $\langle B\rangle$ is not a subgroup of a
conjugate of a vertex group of $\Phi$). Then the visual
decomposition of $(W,S)$ for $\Phi$ (given by the construction for
theorem \ref{MT1}) has an edge group $\langle A\rangle$ such that
$A$ separates $B$ in $\Gamma$ - which is impossible.

For the third part of the theorem, let $\Lambda$ be the
decomposition of $\langle V\rangle$ determined by its action on
$\Phi$. Since $\Lambda$ is compatible with $\Psi$,  if $\Lambda$
is non-trivial, then $M(\Psi,\langle V\rangle)$ is non-empty. But
as $\Psi$ is JSJ, $M(\Psi,\langle V\rangle)$ can only contain
crossing elements. 
\end{proof}

\section{Orbifold groups} \label{Orbi}

If $\Psi$ is the JSJ-decomposition of $(W,S)$ and $\langle
V\rangle$ is a vertex group of $\Psi$ such that $M(\Psi ,\langle
V\rangle)$ contains no crossing members, then $\langle V\rangle$
is indecomposable with respect to its action on the Bass-Serre
tree for any splitting of $W$ over virtually abelian subgroups.
If instead, $M(\Psi,\langle V\rangle)$ contains crossing members,
then we say $\langle V\rangle$ is an {\it orbifold vertex group of
$\Psi$}.

\begin{theorem}\label{orbi}
Suppose $(W,S)$ is a finitely generated Coxeter system, $\Psi$ is
the (visual) JSJ-decomposition of $(W,S)$, and $\langle V\rangle$
is an orbifold vertex group of $\Psi$. Then $\langle V\rangle$
decomposes as $\langle T\rangle\times \langle M\rangle$ where
$T\cup M=V$, $M$ generates a virtually abelian group and the
presentation diagram of $T$ is either a loop of length $\geq 4$
(in which case $T$ generates a group that is virtually a closed
surface group) or the presentation diagram of $T$ is a disjoint
union of vertices and simple paths (in which case $T$ generates a
virtually free group with graph of groups decomposition such that
each vertex group is either $\mathbb Z_2$ or finite dihedral and
each edge group is either trivial or $\mathbb Z_2$).
\end{theorem}

\begin{proof}
As $\Psi$ is JSJ,  the only visual
subgroups of $M(\Psi,\langle V\rangle)$ are crossing. Say
$\langle A\rangle$ and $\langle B\rangle$ are minimal rank
elements of $M(\Psi,\langle V\rangle)$ that cross one another. If
$C\subset V$ such that $\langle C\rangle$ is virtually abelian,
and $C$ separates $V$ in $\Gamma$ then $Rank(C)\geq Rank(A)$. If
additionally $Rank(C)=Rank(A)$ then $\langle C\rangle \in
M(\Psi,\langle V\rangle)$ and so $\langle C\rangle$ is crossing
in $M(\Psi,\langle V\rangle)$.

By lemma \ref{L10}, $A=\{a_1,a_2\}\cup M$ and $B=\{b_1,b_2\}\cup
M$ where $\langle M\rangle$ is virtually abelian and commutes with
$\{a_1,a_2,b_1,b_2\}$, $A$ separates $b_1$ and $b_2$ in $\Gamma$
and $B$ separates $a_1$ and $a_2$ in $\Gamma$. By Proposition
\ref{P11}, $\Gamma -A$ has exactly two components which intersect
$V$ non-trivially, $C_{b_1}$ containing $b_1$ and $C_{b_2}$
containing $b_2$. Then $V -M\subset \{a_1\}\cup C_{b_1}\cup
\{a_2\}\cup C_{b_2}$. By lemma \ref{L8} there is an edge
connecting (cyclically) adjacent members of this union, and no
edge connecting non-adjacent members. In particular, $C_{b_1}\cup
\{a_1,a_2\}$ is a connected subset of $\Gamma$. We show:

\medskip

\noindent (i) $(C_{b_1}\cup \{a_1,a_2\})-\{b_1\}$ has exactly two
components which intersect $V$ non-trivially, one containing $a_1$
and the other containing $a_2$.

\medskip

A path connecting $a_1$ and $a_2$ in $C_{b_1}\cup \{a_1,a_2\}$
avoiding $b_1$ also avoids $b_2$ and $M$, but this is impossible
as $B$ separates $a_1$ and $a_2$ in $\Gamma$. Hence $b_1$
separates $a_1$ and $a_2$ in $(C_{b_1}\cup \{a_1,a_2\})$. Let
$C_{b_1,a_1}$ (respectively $C_{b_1,a_2}$) be the component of
$(C_{b_1}\cup \{a_1,a_2\})-\{b_1\}$ containing $a_1$
(respectively $a_2$). If $z\in V$ is an element of a third
component of $(C_{b_1}\cup \{a_1,a_2\})-\{b_1\}$, then $z\in
\Gamma -B$ and by proposition \ref{P11}, there is a (shortest) path
$\alpha$ in $\Gamma -B$ from $z$ to either $a_1$ or $a_2$. Without
loss, say $\alpha$ connects $z$ to $a_1$. Then $a_2$ is not a
point of $\alpha$. Note that $\alpha$ is not contained in
$C_{b_1}\cup \{a_1,a_2\}$ since $b_1$ separates $z$ and $a_1$ in
$C_{b_1}\cup \{a_1,a_2\}$ and $\alpha$ avoids $b_1$. If $\beta$
is the longest initial segment of $\alpha$ in $C_{b_1}\cup
\{a_1,a_2\}$ with end point $t$, and $s$ follows $t$ on $\alpha$,
then $\beta$ avoids $A$ and $s\not\in M\cup \{a_1,a_2\}\cup
C_{b_1}$. Then $s$ is connected to $z\in V$ by a path in $\Gamma
-A$, implying $s\in C_{b_1}\cup C_{b_2}$. So $s\in C_{b_2}$. But
then the edge $[ts]$ connects points of $C_{b_1}$ and $C_{b_2}$,
which is impossible (recall $A$ separates $C_{b_1}$ and $C_{b_2}$
in $\Gamma$). So (i) is proved.

Observe that $V$ is contained in the union $C_{b_1,a_1}\cup
C_{b_1,a_2}\cup C_{b_2,a_1}\cup C_{b_2,a_2}\cup M\cup
\{b_1,b_2\}$. The only overlap among these sets is
$C_{b_1,a_i}\cap C_{b_2,a_i}=\{a_i\}$. Next we show:

\medskip
The set $\{b_1,a_1\}\cup M $ separates $C_{b_1,a_1}-\{a_1\}$ from $C_{b_2}\cup C_{b_1,a_2}$ in $\Gamma$. By the symmetry of the situation (see figure \ref{jsjfig4}), the same argument implies:

\medskip

\noindent (ii) For $u\in \{b_1,b_2\}$ and $v\in \{a_1,a_2\}$,
$\{u,v\}\cup M$
separates $C_{u,v}-\{v\}$ from $
C_t\cup C_{u,s}$ in $\Gamma$, where $\{t\}=\{b_1,b_2\}-\{u\}$ and
$\{s\}=\{a_1,a_2\}-\{v\}$.

\medskip

\begin{figure}[ht]
\begin{center}
\epsfig{file=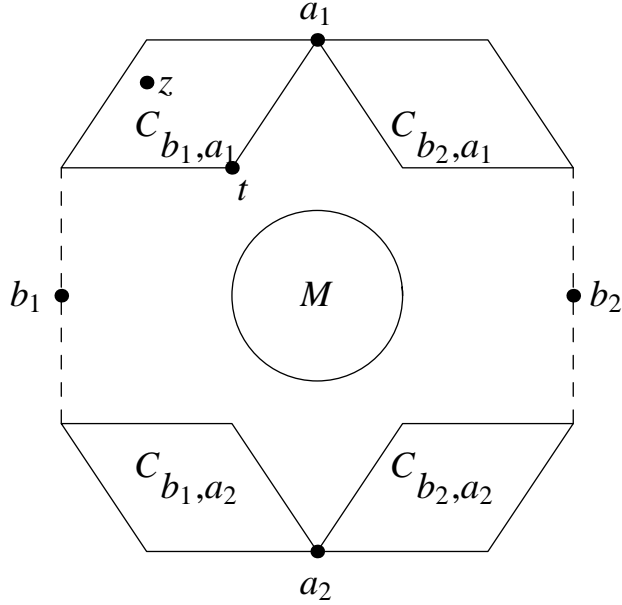}
\caption{Circular decomposition of $V$ in $\Gamma(W,S)$. \label{jsjfig4}}
\end{center}
\end{figure}

\medskip

By lemma \ref{L8} there is an edge from $a_2$
to $C_{b_2}$ so $C_{b_2}\cup C_{b_1,a_2}$ is a connected subset of
$\Gamma -(\{b_1,a_1\}\cup M)$. Suppose $z\in
C_{b_1,a_1}-\{a_1\}$ and $\alpha$ is a path from $z$ to $a_2$ in
$\Gamma -(\{a_1,b_1\}\cup M)$. Let $t$ be the last vertex of
$\beta \equiv$ the longest initial segment of $\alpha$ in
$C_{b_1,a_1}-\{a_1\}$ (then $t\ne a_2$) and $s$ the vertex of
$\alpha$ following $t$. Note that $s\not \in \{a_1,b_1\}\cup M$
by the definition of $\alpha$, $s\not\in C_{b_1,a_1}-\{a_1\}$ by
the definition of $\beta$, and $s$ is not an element of a
component of $(C_{b_1}\cup \{a_1,a_2\})-\{b_1\}$ other than
$C_{b_1,a_1}$, since $b_1$ separates $C_{b_1,a_1}$ from those
components in $C_{b_1}\cup \{a_1,a_2\}$. Hence $s\not \in
C_{b_1}\cup \{a_1,a_2\}\cup M$. The path $\beta$ followed by the edge from $t$ to $s$ begins in $C_{b_1}$ and avoids $A$. This implies $s\in C_{b_1}$, which is nonsense. So, (ii) is proved.


For $u,v$ as in (ii), if $V\cap C_{u,v}\ne \{v\}$, then
$\{u,v\}\cup M\equiv C$ separates $V$ in $\Gamma$. There is no
edge connecting $u$ and $v$ since otherwise, $Rank (\langle
C\rangle)<Rank \langle A\rangle$. Hence $Rank(\langle
C\rangle)=Rank (\langle A\rangle)$, and $\langle C\rangle$ is
crossing in $M(\Psi,\langle V\rangle)$. Say $\langle C\rangle$
crosses $\langle D\rangle$. If $D$ does not separate $u$ and $v$
in $\Gamma$ then $D$ must separate some other unrelated pair,
$m_1$ and $m_2$ in $M$. But then $D$ separates $A$ in $\Gamma$ and $D$ separates $B$ in $\Gamma$.
By lemma \ref{L10}, $D-\{d_1,d_2\}=A-\{m_1,m_2\}$ and $D-\{d_1',d_2'\}=B-\{m_1,m_2\}$. Hence $\{a_1,a_2,b_1,b_2\}\subset D$. As $D$ is virtually abelian and the pairs $a_1,a_2$ and $b_1,b_2$ are unrelated pairs, there is an edge (labeled 2) in $\Gamma$ connecting $a_i$ to $b_j$ for all $i,j\in \{1,2\}$.
This is nonsense since we have assumed $u\in \{b_1,b_2\}$ and $v\in \{a_1,a_2\}$ are unrelated. Instead $D$
separates $u$ and $v$ in $\Gamma$.

If $V\cap C_{b_1,a_1}\ne \{a_1\}$, then proposition \ref{P11}
implies $\Gamma -(\{a_1,b_1\}\cup M)$ has exactly 2 components
which intersect $V$ non-trivially. Now, $\{b_1,a_1\}\cup M$
separates $C_{b_1,a_1}-\{a_1\}$ from (the connected set)
$C_{b_2}\cup C_{b_1,a_2}$ in $\Gamma$. Hence we have:

\medskip

\noindent (iii) If $V\cap C_{b_1,a_1}\ne\{a_1\}$ then
$V\cap (C_{b_1,a_1}-\{a_1\})$ is contained in one of the two components of
$\Gamma -(\{a_1,b_1\}\cup M)$ which intersect $V$ non-trivially (call it $Q_1$)
and $C_{b_2}\cup C_{b_1,a_2}$ is contained in the other (call it $Q_2$).

\medskip

Next we improve (iii) by proving:

\medskip

\noindent (iv) If $V\cap C_{b_1,a_1}\ne\{a_1\}$ then $Q_1$,
the component of $\Gamma -(\{a_1,b_1\}\cup M)$ containing $V\cap (C_{b_1,a_1}-\{a_1\})$, is a component of $C_{b_1,a_1}-\{a_1\}$. (In particular, 
the set $V\cap (C_{b_1,a_1}-\{a_1\})$ is contained in a single component of $C_{b_1,a_1}-\{a_1\}$.)

\medskip

If $K_1$ is a component of $\Gamma -A$ other than $C_{b_1}$ or
$C_{b_2}$ and $\alpha$ is a path connecting $K_1$ to
$C_{b_1,a_1}-\{a_1\}(\subset C_{b_1})$ in $\Gamma
-(\{a_1,b_1\}\cup M)$, then $\alpha$ contains $a_2\in
C_{b_1,a_2}$ which is impossible by (iii). Hence $K_1\cap Q_1=\emptyset$. This implies $Q_1\subset A\cup C_{b_1}\cup C_{b_2}$. But $C_{b_1,a_2}\cup C_{b_2}\subset Q_2$ and $(\{a_1\}\cup M)\cap Q_1=\emptyset$, hence $Q_1\subset C_{b_1}-(\{b_1\}\cup C_{b_1,a_2} )$. Suppose $K_2$ is a
component of $(C_{b_1}\cup \{a_1,a_2\})-\{b_1\}$ other than $C_{b_1,a_1}$ or $C_{b_1,a_2}$ and
$\alpha$ is a path from $K_2$ to $V\cap (C_{b_1,a_1}-\{a_1\})$ ($\subset Q_1$) avoiding
$\{a_1,b_1\}\cup M$. If $\alpha$ leaves $C_{b_1}$, then $\alpha$
contains $a_2$ which is impossible by (iii). Then $\alpha$ is
a path in $C_{b_1}$ passing through $b_1$, which is also
impossible. So, $Q_1\subset C_{b_1,a_1}-\{a_1\}\subset \Gamma-(\{a_1.b_1\}\cup M)$, and (iv) is verified.

If $V\cap C_{b_i,a_j}\ne \{a_j\}$, define $\hat C_{b_i,a_j}$ to be the component of $C_{b_i,a_j}-\{a_j\}$ (equivalently by (iv), the component of $\Gamma-(\{a_j,b_i\}\cup M)$) containing $V\cap(C_{b_i,a_j}- \{a_j\})$.
As $\hat C_{b_i,a_j}\cap(A\cup B)=\emptyset$, $\hat C_{b_i,a_j}$ is a component of $\Gamma- (A\cup B)$. If $V\cap C_{b_i,a_j}=\{a_j\}$, then define $\hat C_{b_i,a_j}\equiv \emptyset$. 

As $C_{b_i,a_j}$ is connected, $a_j$ is connected by an edge to each component of $C_{b_i,a_j}-\{a_j\}$. 
Each component of $C_{b_i,a_j}-\{a_j\}$ is also a component of $C_{b_i}-\{b_i\}$ 
and so is connected by an edge to $b_j$. In particular, if $\hat C_{b_i,a_j}\ne \emptyset$ then it is connected to both $a_j$ and $b_i$ by edges. 


We conclude that $V -(A\cup B)$ is contained in the following
disjoint union of sets, each of which is either trivial, or has nontrivial 
intersection with $V$ and is a component of $\Gamma -(A\cup B)$.
$$\hat C_{b_1,a_1}\cup \hat C_{b_1,a_2}\cup \hat C_{b_2,a_2}\cup
\hat C_{b_2,a_1}$$

If $\hat C_{b_i,a_j}=\emptyset$, then there may be an
edge connecting $b_i$ and $a_j$, but if no such set is empty,
then two members of the following collection are connected by an
edge iff they are cyclically adjacent. I.e. they form a ``loop"
of sets in $\Gamma$:
$$\{a_1\}, \hat C_{b_1,a_1}, \{b_1\}, \hat C_{b_1,a_2},
\{a_2\}, \hat C_{b_2,a_2}, \{b_2\}, \hat C_{b_2,a_1}$$

By the symmetry of the above results, we could define
$C_{a_1,b_1}$. If $V\cap C_{a_1,b_1}\ne \{b_1\}$ then
$\hat C_{a_1,b_1}=\hat C_{b_1,a_1}$ and is the only
component of $\Gamma -(A\cup B)$ containing points of $V$ and
connected to each of $a_1$ and $b_1$ by edges. Similarly for
the other such sets. Remove $\hat C_{u,v}$ from the above loop
if $\hat C_{u,v}=\emptyset$ and then sets are connected by an edge
only if they are cyclically adjacent.

If $\{u,v\}\in \{\{a_1,b_1\}, \{b_1,a_2\}, \{a_2,b_2\},
\{b_2,a_1\}\}$ and $V\cap C_{u,v}\ne\{v\}$, then $\langle
\{u,v\}\cup M\rangle$ crosses say $\langle \{d_1,d_2\}\cup
M\rangle$ in $M(\Psi,\langle V\rangle)$. By (iii) and proposition
\ref{P11}, for some $i$, $d_i\in C_{u,v}-\{v\}$.  By the above
analysis, $d_i$ separates $u$ and $v$ in $C_{u,v}\cup \{u\}$.
Form the sets $C_{d_i,u}$ and $C_{d_i,v}$ as before. Note that
$d_i$ commutes with $M$.  Continue until all remaining $C_{u,v}$
are such that $V\cap C_{u,v}=\{v\}$. The final collection of
sets are all of the singletons of $V-M$ and each commutes with
$M$. Then $T(\equiv V-M)$ and $M$ satisfies the conclusion of our
theorem.

If the presentation diagram for $T$ is a loop of length $\geq 4$,
then theorem 7.16.2 of Beardon's book \cite{Beardon} guarantees
the existence of an $n$-sided hyperbolic polygon whose vertex
angles are $\frac \pi {m_i}$ (in cyclic order) for $m_i$ the edge
labels of the presentation diagram (in cyclic order). Theorem
7.1.4 of \cite{Rat} concludes that the reflection group in this
Coxeter polygon is a Coxeter group with cyclic presentation
diagram and edge labels $m_i$ (in cyclic order). Selberg's lemma
implies this Coxeter group has a torsion free subgroup of finite
index and so $\langle T\rangle$ has a closed surface subgroup of
finite index.

If $T$ is not a loop, a visual decomposition of $\langle
T\rangle$ produces the desired graph of groups decomposition of
$\langle T\rangle$. 
\end{proof}

\begin{figure}[ht]
\begin{center}
\epsfig{file=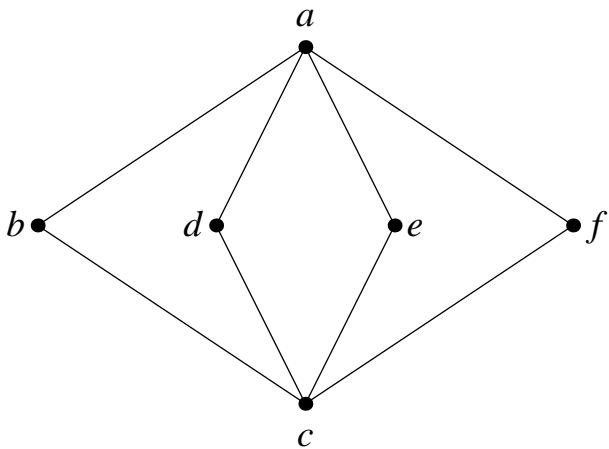}
\caption{$\Gamma(W,S)$ \label{jsjfig5}}
\end{center}
\end{figure}

Note that for a Coxeter system $(W,S)$ $A\subset S$ may separate
$\Gamma(W,S)$ and generate a visual virtually abelian group, but
$A$ may not be a subset of a vertex of a JSJ-decomposition with
virtually abelian edge groups for $W$.

\begin{example} \label{E4}
Let $(W,S)$ be the Coxeter system with presentation diagram shown in figure \ref{jsjfig5} (with all edge labels equal to 2).
Then $W$ has virtually abelian JSJ-decomposition:
$$ \langle a,b,c\rangle\ast _{\langle a,c\rangle} \langle
a,c,d\rangle_{\langle a,c\rangle} \langle a,c,e\rangle _{\langle
a,c\rangle} \langle a,c,f\rangle$$
The set $\{ a,b,c,d\}$ separates $\Gamma$, and generates a
$rank$-2 virtually abelian subgroup of $W$. Each vertex group
of the JSJ-decomposition is $rank$-1 and so cannot contain a conjugate of $\langle a,b,c,d\rangle$.
\end{example}

The group of example 3 (figure \ref{jsjfig3}) has orbifold vertex group $\langle x,y,u,v\rangle$, a virtually free group with presentation diagram two disjoint edges. This group is isomorphic to $\langle x,u\rangle\ast \langle y,v\rangle$, the free product of two (finite) dihedral groups. 

\begin{figure}[ht]
\begin{center}
\epsfig{file=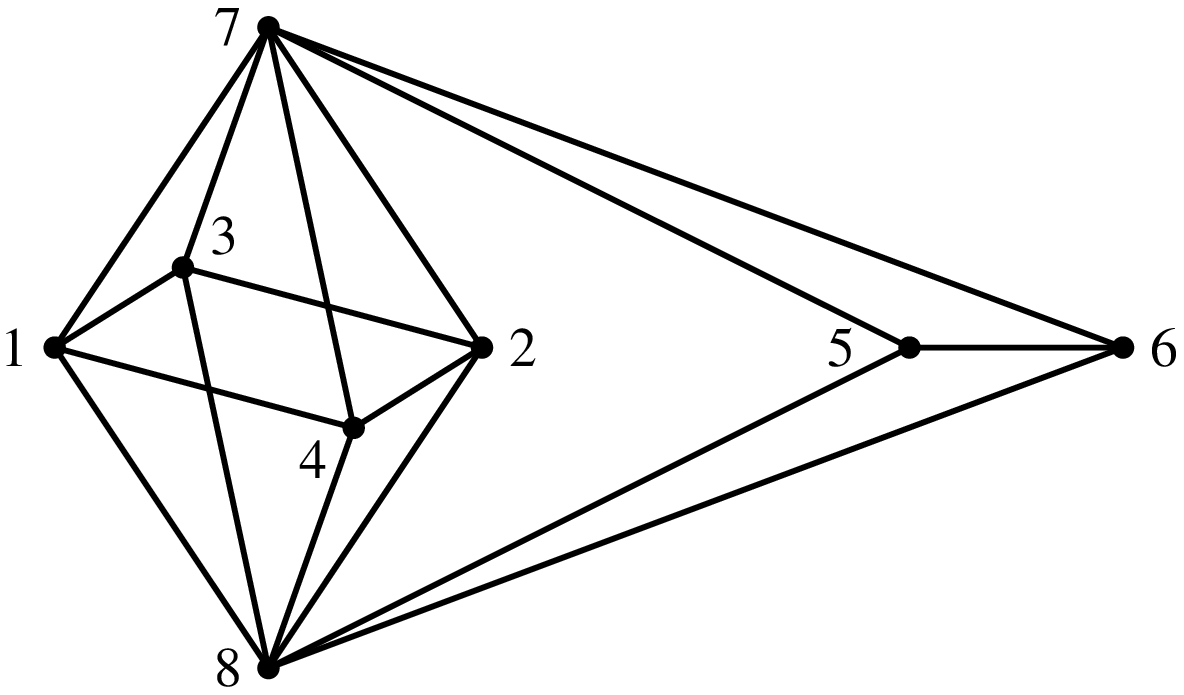}
\caption{$\left(\langle7,8\rangle\times\langle 1,2,3,4\rangle \right) \ast_{\langle 7,8\rangle}\langle 5,6,7,8\rangle$ \label{jsjfig6}}
\end{center}
\end{figure}

\begin{example}\label{E5} 
Let $(W,S)$ be the Coxeter system with presentation diagram shown in figure \ref{jsjfig6} (with all edge labels connecting $\{7,8\}$ to $\{1,2,3,4\}$ equal to 2 and other labels arbitary).
Then $W$ has virtually abelian JSJ-decomposition: 
$$\left(\langle7,8\rangle\times\langle 1,2,3,4\rangle \right) \ast_{\langle 7,8\rangle}\langle 5,6,7,8\rangle$$
In the orbifold vertex group $ \langle7,8\rangle\times\langle 1,2,3,4\rangle$ the virtually abelian splitting subgroups $\langle 1,2\rangle \times \langle 7,8\rangle$ and $\langle 3,4\rangle \times \langle 7,8\rangle$ are crossing. The loop determined by vertices in $\{1,2,3,4\}$ generates a virtual closed surface group. 
\end{example}

If in the preceding example the loop determined by the vertices $\{1,2,3,4\}$ is replaced by a loop of  length $\geq 4$ with arbitrary edge labels and one considers the direct product of such a  Coxeter group with an arbitrary virtually abelian Coxeter group then the resulting Coxeter groups have JSJ decompositions with  orbifold groups that  (non-trivially) realizes all possible orbifold groups of the type as described in part 3) of theorem 1.

\begin {example}
Suppose $(W',T)$ has presentation diagram that is a union of isolated vertices and simple paths. Our goal is to realize $\langle T\rangle$ as an orbifold group in the JSJ decomposition of a 1-ended Coxeter group. Decompose $T$ as $\cup _{i=0}^nT_i$ where the $T_i$ are disjoint and each $T_i$ is the set of vertices of a component of the presentation diagram $\Gamma(W',T)$. If $T_i$ is a singleton $z_i$, then let $a_{i}=b_{i}=z_i$. Otherwise, let $a_{i}$ and $b_{i}$ be the end points of the simple path of $T_i$. 
Let $S=T\cup (\cup_{i=0}^n C_i)$ where (mod $n+1$) $C_i=\{b_{i}, a_{i+1},x_{i},y_{i}\}$ for $i\in \{0,\ldots ,n\}$ and $(W_i, C_i)$ is the Coxeter system with  group $W_i $ isomorphic to $\mathbb D_2\times \mathbb Z_2\times \mathbb Z_2$ (where $\{x_i\}$ generates one of  the $\mathbb Z_2$, $\{y_i\}$ generates the other, and $\{b_i,a_{i+1}\}$ generates the infinite dihedral group $\mathbb D_2$). 

The desired Coxeter system $(W,S)$ has the relations of $(W',T)$ and $(W_i,C_i)$. The JSJ decomposition of $(W,S)$ has irreducible vertex groups $\langle C_i\rangle$ and orbifold vertex group $\langle T\rangle$. The tree for this decomposition has a vertex $V$ (with vertex group $W'\equiv \langle T\rangle$), $n+1$ edges $E_i$, each incident to $V$  (the edge group of $E_i$ is the infinite dihedral group $\langle b_i,a_{i+1} \rangle$), and  $V_i$, the vertex of $E_i$ opposite $V$, has edge group $W_i\equiv \langle C_i\rangle$). See figure 7.

The vertices of $T$ can now be ordered as $t_0,\ldots t_m$ respecting the ordering of the $T_i$ and the internal ordering of the simple paths with $a_i$ preceding $b_i$. If $t_i$ and $t_j$ are not adjacent (mod $m+1$) then $\{ t_i,t_j\}$ separates $T$ in $\Gamma(W,S)$ and defines a crossing splitter. The group $W$ is 1-ended since no subset $A$ of $S$ both separates $\Gamma(W,S)$ and generates a finite subgroup of $W$. 
\end{example}
\begin{figure}[ht]
\begin{center}
\epsfig{file=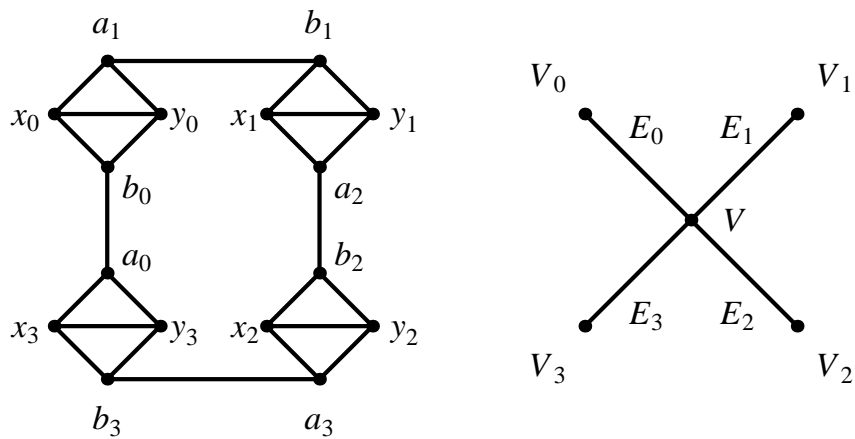}
\caption{Virtually free orbifold group \label{mikefig0619096}}
\end{center}
\end{figure}





\begin{thebibliography}{99}
\bibitem{Beardon} A. F. Beardon, {\it The Geometry of Discrete Groups}, Graduate
Texts in Math., vol. {\bf 91}, Springer Verlag, Berlin,
Heidelberg, and New York, 1983.
\bibitem{Bourbaki} N. Bourbaki, {\it Groupes et Algebres de Lie}
Chapters IV-VI, Hermann, Paris, 1968.
\bibitem{Cox} H. M. S. Coxeter {\it Regular Polytopes}, Dover,
New York, 1973.
\bibitem{Deo} V. V. Deodhar, {\it On the root system of a Coxeter
group}, Commun. Algebra {\bf 10} (1982), 611-630.
\bibitem{DicksDunwoody} W. Dicks and M. J. Dunwoody, {\it Groups
Acting on Graphs}, Cambridge University Press, 1989.
\bibitem{DS} M. J. Dunwoody and M. E. Sageev, {\it JSJ-Splittings
for Finitely Presented Groups over Slender Groups}, Invent. Math.
{\bf 135} (1999), 25-44.
\bibitem{Fo} M. Forester, {\it Deformation and rigidity of simplicial group actions on trees}, Geom. $\&$ Topol. {\bf 6} (2002), 219-267.
\bibitem{FP} K. Fujiwara and P. Papasoglu, {\it JSJ-decompositions of finitely presented groups and complexes of groups}, Geom. Funct. Anal. {\bf 16} (2006) no.1, 70-125.
\bibitem{GL} V. Guirardel and G. Levitt, {\it A general construction of JSJ decompositions}, Geometric Group Theory, Trends in Mathematics, 65-73, 2007 Birkh\"auser Verlag, Basel, Switzerland.
\bibitem{JS} W. Jaco and P. Shalen, {\it Seifert fibered spaces in 3-manifilds}, Memoirs of the AMS {\bf 220}, 1979. 
\bibitem{J} K. Johannson, {\it Homotopy equivalences of 3-manifolds with boundary}, LNM Springer Verlag, 1978. 
\bibitem{Krammer} D. Krammer, {\it The conjugacy problem for
Coxeter groups}, Thesis, Universiteit Utrecht, Netherlands, 1994.
\bibitem{MT} M. L. Mihalik and S. Tschantz, {\it
Visual Decompositions of Coxeter groups}, Groups Geom.  Dyn. {\bf 3} (2009), 173-198. .
\bibitem{Rat} J. G. Ratcliffe, {\it Foundations of Hyperbolic
Manifolds}, Graduate Texts in Math., vol. {\bf 149},
Springer-Verlag, Berlin, Heidelberg, and New York, 1994.
\bibitem{RS} E. Rips and Z. Sela, {\it Cyclic Splittings of
Finitely Presented Groups and the Canonical JSJ-Decomposition},
Ann. of Math (2) {\bf 146} (1997) no. 1, 53-109.
\bibitem{SS} P. Scott and G. A. Swarup, {\it Regular neighbourhoods and canonical decompositions for groups}, Ast\'erisque {\bf 289} (2003).
\bibitem{S} Z. Sela, {\it Structure and rigidity in (Gromov) hyperbolic groups and discrete groups in rank 1 Lie groups II}, Geom. Funct. Anal. {\bf 7}, (1997) no. 3, 561-594.  
\bibitem{Serre} J. P. Serre, {\it Trees}, Springer-Verlag, New
York, 1980.
\bibitem{Wald} F. Waldhausen, {\it On the determination of some bounded 3-manifolds by their fundamental groups alone}, Proc. Sympos. on Topology and its Applications, 331-332, Beograd.
\end{thebibliography}
\end{document}